\documentclass[11pt]{article}
\usepackage{amscd}
\usepackage{amsfonts}
\usepackage{amsmath}
\usepackage{amssymb}
\usepackage{amsthm}
\usepackage{bbm}
\usepackage{CJK}
\usepackage{fancyhdr}
\usepackage{graphicx}
\usepackage{hyperref}
\usepackage{indentfirst}
\usepackage{latexsym}
\usepackage{mathrsfs}
\usepackage{xypic}
\usepackage{amsmath,amssymb,amsfonts,amsthm,fancyhdr}
\usepackage{epsfig,graphicx,picins,picinpar,subfigure}
\usepackage{pstricks}
\usepackage{amssymb}
\usepackage{xypic}
\usepackage[compress]{cite}

\newtheorem{theorem}{Theorem}[section]
\newtheorem{prop}[theorem]{Proposition}
\newtheorem{lemma}[theorem]{Lemma}
\newtheorem{coro}[theorem]{Corollary}
\newtheorem{prop-def}{Proposition-Definition}[section]

\theoremstyle{definition}
\newtheorem{defn}[theorem]{Definition}

\newtheorem{remark}[theorem]{Remark}
\newtheorem{exam}[theorem]{Example}

\usepackage[top=1in,bottom=1in,left=1.25in,right=1.25in]{geometry}
\def\<{\langle}
\def\>{\rangle}

\date{\today}
\begin{document}
\renewcommand{\baselinestretch}{1.2}
\renewcommand{\arraystretch}{1.0}
\title{\bf  Twisted Rota-Baxter family operators on  Hom-associative algebras }
\author{{\bf Wen Teng$^{1,*}$, Yunpeng Xiao$^{2}$}\\
{\small 1   School of Mathematics and Statistics, Guizhou University of Finance and Economics,} \\
{\small  Guiyang  550025, P. R. of China}\\
{\small 2  School of Mathematical Sciences, Guizhou Normal University,  Guiyang  550025, P. R. of China}\\
  {\small   * Corresponding author: tengwen@mail.gufe.edu.cn}\\
 }

 \maketitle
\begin{center}
\begin{minipage}{13.cm}

{\bf Abstract:}
In this paper, we first define twisted Rota-Baxter family operators
  on  Hom-associative algebras indexed by
a semigroup $\Omega$. Then we introduce and study Hom-NS-family algebras as the underlying structures of twisted Rota-Baxter family operators. Meanwhile, We show that a Hom-NS-family algebra induces an ordinary Hom-NS-algebra  on the tensor product with the semigroup algebra.
Moreover,  we define the cohomology of a twisted   Rota-Baxter family operator. This cohomology
can also be viewed as the   cohomology of a
certain Hom-$\Omega$-associative algebra with coefficients in a suitable
bimodule.  Finally, we examine
deformations of twisted  Rota-Baxter family operators and demonstrate that  they are governed
by the aforementioned cohomology. The concept
of Nijenhuis elements  linked to a twisted   Rota-Baxter family operator is
introduced to  provide a sufficient condition for its rigidity.
 \smallskip

 \smallskip

 \end{minipage}
 \end{center}
 \normalsize\vskip0.5cm

\section{Introduction}
\def\theequation{\arabic{section}. \arabic{equation}}
\setcounter{equation} {0}

Rota-Baxter operators, commonly known as Baxter operators, embody the algebraic generalization of the integral operator.
Rota-Baxter operators initially emerged in Baxter's work addressing the problem of probability \cite{Baxter} and were subsequently studied extensively by Rota \cite{Rota}  in relation to combinatorics.
 A significant number of scholarly works have been dedicated to exploring the multifaceted nature of Rota-Baxter operators, delving into the realms of both pure and applied mathematics \cite{Gubarev,Bai12,Ebrahimi-Fard07,Ebrahimi-Fard08}.    For further details on Rota-Baxter operators, see  \cite{Guo12}.
In \cite{Uchino}, Uchino introduced generalized Rota-Baxter operators on bimodules over an associative algebra. These operators are called  twisted Rota-Baxter operators or
 twisted $\mathcal{O}$-operators, which can be seen as a non-commutative analog of  twisted  Poisson structures \cite{Severa}.
Uchino also noted that a   twisted Rota-Baxter operator gives rise to an NS-algebra of Leroux \cite{Leroux} in a manner analogous to how a Rota-Baxter operator induces a dendriform algebra structure.
 See   \cite{Das2022} for further details on NS-algebras.

The Rota-Baxter family algebra, introduced by Guo \cite{Guo2009}, constitutes a further generalization of the Rota-Baxter algebra. In this context, the standard operation is substituted with a family of operations indexed by a semigroup $\Omega$, which arises naturally in the renormalization of quantum field theory.
Recently, various family algebraic structures have been defined, as referenced in \cite{Foissy2021,Foissy2024,Zhang2019,Zhang2020}.
 In particular, Das \cite{Das22} presented the notion of an $\mathcal{O}$-operator family on associative algebras, which serves as  a generalization of the Rota-Baxter family.
There, the author established cohomology and deformation theory of  twisted   $\mathcal{O}$-operator family on associative algebras and NS-family algebras.
Utilizing Das's method,
Liu and Zheng developed the cohomology and deformation theory of  twisted   $\mathcal{O}$-operator families on    Leibniz algebras  in \cite{Liu}.

Hom-type algebras originated from $q$-deformations of   Witt and Virasoro
algebras \cite{Hu}, playing a significant role in physics and conformal field theory.
Hartwig, Larsson, and Silvestrov \cite{Hartwig} introduced the concept of Hom-Lie algebras, which generalize traditional Lie algebras by twisting the usual Jacobi identity with a homomorphism.
The concept of Hom-associative algebras was first introduced in \cite{Makhlouf}.
Within this framework, various classical structures have been generalized. A significant amount of research has been dedicated to Hom-type algebras, partly due to the potential of establishing a general framework that enables the production of numerous types of natural deformations, which are of interest to both mathematicians and physicists.
Recently, the cohomology and deformations of Hom-Lie and Hom-associative algebras were investigated in \cite{Ammar2011,Sheng2012,Hurle2019}. It is important to note that the cohomology and deformations
of $\mathcal{O}$-operators on Hom-Lie and Hom-associative algebras  were addressed in \cite{Mishra2020,Chtioui2022}.

Drawing inspiration from the foundational work on (twisted) $\mathcal{O}$-operator families \cite{Das22}, this paper presents a generalization and explores the connections between pertinent findings on associative algebras and Hom-associative algebras. Consequently, this paper aims  to delve into the realm of twisted Rota-Baxter family operators on Hom-associative algebras, which constitute a broader category encompassing associative algebras.
For this purpose,  we  initially present the notion of   twisted Rota-Baxter family operators on  Hom-associative algebras  indexed by a semigroup $\Omega$, and provide several illustrative examples. Then, we define Hom-NS-family algebras,  which are associated with twisted Rota-Baxter family operators. We demonstrate that a  Hom-NS-family algebra   induces an ordinary Hom-NS-algebras. Furthermore, we construct Hom-NS-family algebras utilizing Rota-Baxter family operators of arbitrary weights,  twisted Rota-Baxter family operators and Nijenhuis family operators on Hom-associative algebras. Finally, the cohomology and deformations of  twisted Rota-Baxter family operators on  Hom-associative algebras are discussed.

The paper's structure is presented as follows.
In Section  \ref{sec:Twisted Rota-Baxter}, we introduce twisted Rota-Baxter family operators on Hom-associative algebras, and give some examples.
In Section \ref{sec:Hom-NS-family algebras}, we introduce Hom-NS-family algebras as the underlying structure
of twisted Rota-Baxter family operators. We demonstrate that a Hom-NS-family algebra gives rise to an ordinary Hom-NS-algebra structure on the tensor product with the semigroup algebra.
 In Section \ref{sec:Cohomology},  we introduce Hom-$\Omega$-associative algebra (called Hom-associative algebra  with a semigroup $\Omega$) and its cohomology. This induces the cohomology of twisted Rota-Baxter family operators, which will be used in Section \ref{sec:Deformations} to study the deformations of twisted Rota-Baxter family operators.

\section{Twisted Rota-Baxter family  operators on  Hom-associative algebras }\label{sec:Twisted Rota-Baxter}
\def\theequation{\arabic{section}.\arabic{equation}}
\setcounter{equation} {0}

In this section, we introduce notion  of twisted Rota-Baxter family  operators on  Hom-associative algebras as family analogues of  twisted Rota-Baxter   operators on  Hom-associative algebras. For this,
 we first recall some basic notions about Hom-associative algebras.
 We work over a field $\mathbb{K}$ of
characteristic zero.

\begin{defn} \cite{Makhlouf}
A Hom-associative algebra is a Hom-vector space $(L,p)$, consisting of a vector space $L$ and a linear map $p$, together with a bilinear operation $\mu:L\times L\rightarrow L, (x,y)\mapsto x\cdot y$, that satisfies
\begin{align*}
&p(x\cdot y)=p(x)\cdot p(y),\\
&p(x)\cdot(y\cdot z)=(x\cdot y)\cdot p(z), ~~\forall x,y,z\in L.
\end{align*}
\end{defn}

\begin{defn}
Let $(L,\mu, p)$ be a Hom-associative algebra and $(V,q)$ a Hom-vector space. Let
 $\mu_l:L\times V\rightarrow V, (x,u)\mapsto x\cdot_l u$ and $ \mu_r:V\times L\rightarrow V, (u,x)\mapsto u\cdot_r x$ be two linear maps.
 We say that the
tuple $(V,\mu_l,\mu_r,q)$ is a bimodule of $(L,\mu, p)$ if for any $x,y\in L$ and $u\in V$,
\begin{align*}
&q(x\cdot_l u)=p(x)\cdot_l q(u),~~~~~q(u\cdot_r x)=q(u)\cdot_r p(x),\\
&q(u)\cdot_r (x\cdot y)=(u\cdot_r x)\cdot_r p(y),~p(x)\cdot_l(u\cdot_r y)=(x\cdot_l u)\cdot_r p(y),~ p(x)\cdot_l(y\cdot_l u)=(x\cdot y)\cdot_l q(u).
\end{align*}
\end{defn}

Let $(L,\mu, p)$  be a Hom-associative algebra and $(V,\mu_l,\mu_r,q)$ be a bimodule. The
cohomology of  $L$   with coefficients in $V$ \cite{Ammar2011,Hurle2019} is the cohomology of the cochain complex $C^n_{HA}(L,V)=\{f:L^{\otimes n}\rightarrow V~|~q\circ f=f\circ p^{\otimes n}\}$, $(n\geq 1)$ with the differential
$\partial_{HA}^n:C^n_{HA}(L,V)\rightarrow C^{n+1}_{HA}(L,V)$
 defined by
\begin{align*}
(\partial_{HA}^nf)(x_1,x_2,\ldots,x_{n+1})=&p^{n-1}(x_1)\cdot_l f(x_2,\ldots,x_{n+1})+(-1)^{n+1} f(x_1,\ldots,x_{n})\cdot_r p^{n-1}(x_{n+1})+\\
&\sum_{i=1}^n(-1)^if(p(x_1),\ldots,p(x_{i-1}),x_i\cdot x_{i+1},p(x_{i+2}),\ldots,p(x_{n+1})),
\end{align*}
for $x_1,x_2,\ldots,x_{n+1}\in L.$  An $f\in C^n_{HA}(L,V)$ is called an $n$-cocycle if $\partial_{HA}^nf=0$ and $g$ is called an $n$-coboundary if there exists an $(n-1)$-cochain $f$ such that $\partial_{HA}^{n-1}f=g$. Denote by
$Z^n_{HA}(L,V)$ and $B^n_{HA}(L,V)$ the subspaces of $n$-cocycles and $n$-coboundaries, respectively. The corresponding cohomology groups are denoted by
$H^n_{HA}(L,V)=\frac{Z^n_{HA}(L,V)}{B^n_{HA}(L,V)}$.

In the following, we introduce the twisted Rota-Baxter family  operators on  Hom-associative algebras. Let
$(L,\mu, p)$  be a Hom-associative algebra and $(V,\mu_l,\mu_r,q)$ be a bimodule over it. Suppose that
$\Phi\in C^2_{HA}(L,V)$  is a 2-cocycle of $L$ with coefficients in $V$, i.e. $\Phi:L\otimes L\rightarrow V$ is a bilinear
map satisfying
\begin{align}
&q \Phi(x_1,x_2,x_3)=\Phi(p(x_1),p(x_2),p(x_3)),\label{2.1}\\
&p^{}(x_1)\cdot_l \Phi(x_2,x_{3})- \Phi(x_1,x_{2})\cdot_r p^{}(x_{3})-\Phi(x_1\cdot x_{2},p(x_{3}))+\Phi(p(x_{1}), x_2\cdot x_{3})=0.\label{2.2}
\end{align}

Let $\Omega$ be a fixed semigroup.

\begin{defn}
A collection $\{R_\alpha:V\rightarrow L\}_{\alpha\in\Omega}$ of linear maps is said to be an
$\Phi$-twisted Rota-Baxter family operator   (on $V$ over the Hom-associative algebras $L$) if $R_\alpha$ satisfies
\begin{align}
&R_\alpha\circ q=p\circ R_\alpha,\label{2.3}\\
&R_\alpha u\cdot R_\beta v=R_{\alpha\beta}(R_\alpha u \cdot_l v+u\cdot_rR_\beta v+\Phi(R_\alpha u, R_\beta v)),\label{2.4}
\end{align}
for $u,v\in V$ and $\alpha,\beta\in \Omega.$
\end{defn}

\begin{defn} \label{defn:2.4}
 Let $(L',\mu', p')$  be another  Hom-associative algebra and $(V',\mu'_l,\mu'_r,q')$ be a bimodule of
$L'$. Suppose that $\Phi'\in C^2_{HA}(L',V')$ is a 2-cocycle and $\{R'_\alpha:V'\rightarrow L'\}_{\alpha\in\Omega}$ is an $\Phi'$-twisted Rota-Baxter family operator.
A morphism of twisted Rota-Baxter family operators from $\{R_\alpha:V\rightarrow L\}_{\alpha\in\Omega}$ to $\{R'_\alpha:V'\rightarrow L'\}_{\alpha\in\Omega}$ consists of a pair $(\phi,\varphi)$ of a Hom-associative algebra morphism
$\varphi:(L,\mu, p)\rightarrow(L',\mu', p')$ and a linear map $\phi:V\rightarrow V'$ satisfying
\begin{align}
&\varphi\circ R_\alpha=R'_\alpha\circ\phi,~~~\phi\circ\Phi=\Phi'\circ(\varphi\otimes\varphi),~~~ \phi\circ q=q'\circ\phi, \label{2.5}\\
&\phi\circ \mu_l=\mu'_l \circ(\varphi\otimes\phi),~~\phi\circ \mu_r=\mu'_r \circ(\phi\otimes\varphi).\label{2.6}
\end{align}
\end{defn}

\begin{exam} \label{exam:2.5}
An $\Phi$-twisted Rota-Baxter family operator is a   relative Rota-Baxter family  operator when
$\Phi=0$.
\end{exam}

\begin{exam}
Any $\Phi$-twisted Rota-Baxter operator   is  an $\Phi$-twisted Rota-Baxter family operator indexed by the trivial semigroup reduced to one element.
\end{exam}

\begin{exam} \label{exam: tensor product}
Let $(L,\mu, p)$  be a Hom-associative algebra.  Then the tensor product $L\otimes \mathbb{K}\Omega$
is a Hom-associative algebra with the multiplication $\bar{\cdot}$ and  structure map $\bar{p}$
\begin{align*}
&(x\otimes \alpha)\bar{\cdot}(y\otimes \beta)=x\cdot y\otimes \alpha\beta ~~~~\text{and}\\
&\bar{p}(x\otimes \alpha)=p(x)\otimes \alpha, ~~\forall x,y\in L, \alpha,\beta\in\Omega.
\end{align*}
Moreover, $L$ is a bimodule of $L\otimes \mathbb{K}\Omega$ with the left and right actions
\begin{align*}
&(x\otimes \alpha)\hat{\cdot}_l y =x\cdot y,~~ y\hat{\cdot}_r(x\otimes \beta)=y\cdot x.
\end{align*}
Then we can prove that the map $\hat{\Phi}:(L\otimes \mathbb{K}\Omega)\otimes(L\otimes \mathbb{K}\Omega)\rightarrow L\otimes \mathbb{K}\Omega$ given by
$\hat{\Phi}(x\otimes \alpha,y\otimes \beta)=-x\cdot y$ is a 2-cocycle in the cohomology of $L\otimes \mathbb{K}\Omega$ with coefficients in
$L$. And the collection $\{Id_\alpha:L\rightarrow  L\otimes \mathbb{K}\Omega\}_{\alpha\in\Omega}$ defined by $Id_\alpha(x)=x\otimes\alpha$ for
$x\in L$ is an $\hat{\Phi}$-twisted Rota-Baxter family operator.
\end{exam}

\begin{exam}  \label{exam:2.8}
Let $(L,\mu, p)$  be a Hom-associative algebra   and $\{N_\alpha:L\rightarrow L\}_{\alpha\in\Omega}$ a Nijenhuis family operator
on it, i.e.
\begin{align*}
& N_\alpha x \cdot N_\beta y=N_{\alpha\beta}( N_\alpha x\cdot y+x\cdot  N_\beta y- N_{\alpha\beta}x\cdot y),\\
&p\circ  N_\alpha= N_\alpha\circ p, ~~\forall x,y\in L, \alpha,\beta\in\Omega.
\end{align*}
In this case, $L\otimes \mathbb{K}\Omega$ carries a new  Hom-associative algebra structure with the multiplication ${\cdot}_N$ and structure map $\bar{p}$,
\begin{align*}
&(x\otimes \alpha){\cdot}_N(y\otimes \beta)=N_\alpha x\cdot y\otimes \alpha\beta+ x\cdot N_\beta y\otimes \alpha\beta-N_{\alpha\beta}x\cdot y,\\
&\bar{p}(x\otimes \alpha)=p(x)\otimes \alpha, ~~\forall x,y\in L, \alpha,\beta\in\Omega.
\end{align*}
Moreover, $L$ is a bimodule of $(L\otimes \mathbb{K}\Omega,{\cdot}_N,\bar{p})$ with the left and right actions
\begin{align*}
&(x\otimes \alpha){\cdot}_l^N y =N_\alpha x\cdot y,~~ y{\cdot}_r^N(x\otimes \beta)=y\cdot N_\alpha x.
\end{align*}
It is easy to prove that the map $\Phi_N:(L\otimes \mathbb{K}\Omega)\otimes(L\otimes \mathbb{K}\Omega)\rightarrow L\otimes \mathbb{K}\Omega$ given by
$\Phi_N(x\otimes \alpha,y\otimes \beta)=-N_{\alpha\beta}(x\cdot y)$ is a 2-cocycle in the cohomology of $L\otimes \mathbb{K}\Omega$ with coefficients in
$L$. And the collection $\{Id_\alpha:L\rightarrow  L\otimes \mathbb{K}\Omega\}_{\alpha\in\Omega}$ defined by $Id_\alpha(x)=x\otimes\alpha$ for
$x\in L$ is an $\Phi_N$-twisted Rota-Baxter family operator.
\end{exam}

Consider the Hom-associative algebra $(L\otimes \mathbb{K}\Omega,\bar{\cdot},\bar{p})$ given in Example \ref{exam: tensor product}, we can prove that
$V\otimes \mathbb{K}\Omega$  is a bimodule of $L\otimes \mathbb{K}\Omega$  with the actions
\begin{align*}
&(x\otimes \alpha)\bar{\cdot}_l (u\otimes \beta) =x\cdot_l u \otimes\alpha\beta,~~(u\otimes \beta) \bar{\cdot}_r(x\otimes \alpha) =u\cdot_r x\otimes \beta \alpha.
\end{align*}
Then the map
$\bar{\Phi}:(L\otimes \mathbb{K}\Omega)\otimes(L\otimes \mathbb{K}\Omega)\rightarrow V\otimes \mathbb{K}\Omega$ given by
$\bar{\Phi}(x\otimes \alpha,y\otimes \beta)=\Phi(x,y)\otimes \alpha\beta$
is a 2-cocycle in the cohomology of $L\otimes \mathbb{K}\Omega$ with coefficients in $V\otimes \mathbb{K}\Omega$. Under
this assumption, we have the following conclusion.

 \begin{prop} \label{prop:2.9}
Let $\{R_\alpha:V\rightarrow L\}_{\alpha\in\Omega}$  be an $\Phi$-twisted Rota-Baxter family operator.
Then the map
\begin{align*}
& \bar{R}:V\otimes \mathbb{K}\Omega\rightarrow L\otimes \mathbb{K}\Omega,~~ \bar{R}(u\otimes\alpha)= R_\alpha u\otimes\alpha
\end{align*}
is an $\bar{\Phi}$-twisted Rota-Baxter operator on $V\otimes \mathbb{K}\Omega$ over the Hom-associative algebra $(L\otimes \mathbb{K}\Omega,\bar{\cdot},\bar{p})$.
 \end{prop}

  \begin{proof}
  For all $u,v\in V$ and $\alpha,\beta\in \Omega,$ we have
  \begin{align*}
 \bar{R}  \bar{q}(u\otimes \alpha)= \bar{R} (q(u)\otimes \alpha)=R_\alpha q(u)\otimes \alpha=p(R_\alpha u)\otimes \alpha=\bar{p}(R_\alpha u\otimes \alpha)=\bar{p}(\bar{R} ( u\otimes \alpha)).
\end{align*}
Similarly, we obtain
\begin{align*}
&\bar{R}  ( u\otimes \alpha)\bar{\cdot} \bar{R}  (v\otimes \beta)=(R_\alpha u\otimes\alpha)\bar{\cdot}(R_\beta v\otimes \beta)\\
&=(R_\alpha u {\cdot} R_\beta v)\otimes \alpha\beta=R_{\alpha\beta}(R_\alpha u \cdot_l v+u\cdot_rR_\beta v+\Phi(R_\alpha u, R_\beta v))\otimes \alpha\beta\\
&=\bar{R} (R_\alpha u \cdot_l v+u\cdot_rR_\beta v+\Phi(R_\alpha u, R_\beta v) \otimes \alpha\beta)\\
&=\bar{R} (R_\alpha u \cdot_l v \otimes \alpha\beta+u\cdot_rR_\beta v \otimes \alpha\beta+\Phi(R_\alpha u, R_\beta v ) \otimes \alpha\beta)\\
&=\bar{R}\big(\bar{R}( u\otimes \alpha) \bar{\cdot}_l (v\otimes \beta)+( u\otimes \alpha)\bar{\cdot}_r\bar{R} (v\otimes \beta)+\bar{\Phi}(\bar{R}  ( u\otimes \alpha), \bar{R }(v\otimes \beta))\big).
\end{align*}
Therefore, $\bar{R}$ is an $\bar{\Phi}$-twisted Rota-Baxter operator.
   \end{proof}

Next, we characterize an $\Phi$-twisted Rota-Baxter family  operator by its graph. Let
$(L,\mu, p)$  be a Hom-associative algebra  and $(V,\mu_l,\mu_r,q)$ be a bimodule.  Suppose that $\Phi$
is a 2-cocycle, then the direct sum $L\oplus V$ carries a Hom-associative algebra structure with
the multiplication $\cdot_\Phi$ amd structure map $p\oplus q$ given by
\begin{align*}
&(x,u)\cdot_\Phi(y,v)=(x\cdot y, x\cdot_lv+u\cdot_ry+\Phi(x,y)),\\
&p\oplus q(x,u)=(p(x),q(u))~~~\forall x,y\in L, u,v\in V.
\end{align*}
This is called the $\Phi$-twisted semidirect product Hom-associative algebra, which is denoted
by $L\ltimes_\Phi V$.

Let $(L,\mu, p)$  be a Hom-associative algebra. A collection $\{H_\alpha\}_{\alpha\in \Omega}$ of subspaces of $L$ is said to
be a Hom-family subalgebra of $L$ if $p(H_\alpha)\subseteq H_\alpha$ and $H_\alpha\cdot H_\beta\subseteq H_{\alpha\beta}$ for $\alpha,\beta\in \Omega$.

 \begin{prop}
A collection $\{R_\alpha:V\rightarrow L\}_{\alpha\in\Omega}$  of linear maps is an $\Phi$-twisted Rota-Baxter family  operator if and only if the collection of graphs
$\{Gr(R_\alpha)\}_{\alpha\in\Omega}$ is a Hom-family subalgebra of the $\Phi$-twisted semidirect product Hom-associative algebra $L\ltimes_\Phi V$.
 \end{prop}

 \begin{proof}
Let $(R_\alpha u,u)\in Gr(R_\alpha)$, then
\begin{align*}
&p\oplus q(R_\alpha u,u)=(pR_\alpha u,q(u)).
\end{align*}
On the other hand, for any $(R_\alpha u,u)\in Gr(R_\alpha),(R_\beta v,v)\in Gr(R_\beta)$, we have
\begin{align*}
&(R_\alpha u,u)\cdot_\Phi(R_\beta v,v)=(R_\alpha u\cdot R_\beta v, R_\alpha u\cdot_lv+u\cdot_r R_\beta v+\Phi(R_\alpha u,R_\beta v)).
\end{align*}
Then the collection $\{R_\alpha:V\rightarrow L\}_{\alpha\in\Omega}$   is an $\Phi$-twisted Rota-Baxter family operator if and only if $Gr(R_{\alpha\beta})$ is a Hom--family subalgebra of  $L\ltimes_\Phi V$.
 \end{proof}

\section{Hom-NS-family algebras} \label{sec:Hom-NS-family algebras}
\def\theequation{\arabic{section}.\arabic{equation}}
\setcounter{equation} {0}

In this section, we introduce Hom-NS-family algebras as the family analogue of Hom-NS-algebra. This notion of Hom-NS-family algebra generalizes Hom-dendriform family algebras. We
observed that a   Hom-NS-family algebra induces an ordinary Hom-NS-algebra. Finally, we find the
relations between twisted  Rota-Baxter family operators and Hom-NS-family algebras.

\begin{defn}
 A Hom-NS-algebra is a    quintuple $(G,\prec,\succ, \curlyvee,p)$ consists of a vector space
$G$ together with three bilinear operations $\prec,\succ, \curlyvee:G\otimes G\rightarrow G$  and a linear map $p: G\rightarrow G$   satisfying for $x,y,z\in G$,
 \begin{align}
&p(x\prec y)=p(x)\prec p(y),~~~p(x\succ y)=p(x)\succ p(y),~~~p(x\curlyvee y)=p(x)\curlyvee p(y),\label{3.1}\\
&(x\prec y) \prec p(z)=p(x)\prec (y\prec z+y\succ z+y\curlyvee z),\label{3.2}\\
&(x\succ y)\prec p(z)=p(x)\succ (y\prec z),\label{3.3}\\
&(x\prec y+x\succ y+x\curlyvee y)\succ p(z)=p(x)\succ(y\succ z),\label{3.4}\\
&(x\prec y+x\succ y+x\curlyvee y)\curlyvee p(z)+(x\curlyvee y)\prec p(z)\nonumber\\
&~~~~~~~~~~~~~~~~~~~~~~~~~~~~~~~~~~~=p(x)\succ(y\curlyvee z)+p(x)\curlyvee (y\prec z+y\succ z+y\curlyvee z).\label{3.5}
\end{align}
\end{defn}

\begin{remark}
(i) In \cite{Liu2023}, the  authors studied   BiHom-NS-algebra $(G,\prec,\succ, \curlyvee,p,q)$ , and we can get the concept of Hom-NS-algebra if the structural maps of $G$ is $p=q$.

(ii)  Let   $(G,\prec,\succ, \curlyvee,p)$  be a  Hom-NS-algebra. Then $(G,\ast,p)$ is a Hom-associative algebra, where
 $ x\ast y=x\prec y+x\succ y+x\curlyvee y,~~\forall x,y\in G.$
\end{remark}

\begin{exam}
(i)  Any  NS-algebra is a   Hom-NS-algebra with  $p=Id$.

(ii) Any Hom-dendriform algebra   is a Hom-NS-algebra with $\curlyvee=0$.

(iii)
 Let $(G,\prec,\succ, \curlyvee)$   be an  NS-algebra    and $p:G\rightarrow G$ be an NS-algebra morphism. Define
 $\prec_p,\succ_p, \curlyvee_p:G\otimes G\rightarrow G$ by
$$x\prec_p y=p(x)\prec p(y),~~~x\succ_p y=p(x)\succ p(y),~~~x\curlyvee_p y=p(x)\curlyvee p(y),$$
for all $x,y\in G$. Then $(G,\prec_p,\succ_p, \curlyvee_p,p)$   is a  Hom-NS-algebra, called the
Yau twist of $(G,\prec,\succ, \curlyvee)$.
\end{exam}

 The family version of Hom-NS-algebra is given by the following.

\begin{defn}
  A  Hom-NS-family algebra is  a  linear space $G$ together with
a family of bilinear operations  $\{\prec_\alpha,\succ_\alpha,\curlyvee_{\alpha,\beta}:G\times G\rightarrow G\}_{\alpha,\beta\in \Omega}$
and a linear map  $p:G\rightarrow G$ satisfying the following
conditions
 \begin{align}
&p(x\prec_\alpha y)= p(x)\prec_\alpha p(y),~~p(x\succ_\alpha y)= p(x)\succ_\alpha p(y),~~p(x\curlyvee_{\alpha,\beta} y)= p(x)\curlyvee_{\alpha,\beta} p(y),\label{3.6}\\
& p(x)\prec_{\alpha\beta}(y\prec_\beta z+y\succ_\alpha z+y\curlyvee_{\alpha,\beta}z)=(x\prec_\alpha y)\prec_\beta p(z),\label{3.7}\\
&(x\succ_\alpha y)\prec_\beta p(z)=p(x)\succ_\alpha(y\prec_\beta z),\label{3.8}\\
&(x\prec_\beta y+x\succ_\alpha y+x\curlyvee_{\alpha,\beta}y)\succ_{\alpha\beta}p(z)=p(x)\succ_\alpha(y\succ_\beta z),\label{3.9}\\
&p(x)\succ_\alpha(y\curlyvee_{\beta,\gamma}z)+p(x)\curlyvee_{\alpha,\beta\gamma}(y\prec_\gamma z+y\succ_\beta z+y\curlyvee_{\beta,\gamma}z)\nonumber\\
&=(x\curlyvee_{\alpha,\beta}y)\prec_\gamma p(z)+(x\prec_\beta y+x\succ_\alpha y+x\curlyvee_{\alpha,\beta}y)\curlyvee_{\alpha\beta,\gamma} p(z), \label{3.10}
\end{align}
for all $x,y,z\in G$  and  $\alpha,\beta,\gamma\in \Omega.$

A   Hom-NS-family algebra  is called a regular Hom-NS-family algebra if
$p$ is bijective map.
\end{defn}

A morphism $f: (G,\{\prec_\alpha,\succ_\alpha,\curlyvee_{\alpha,\beta}\}_{\alpha,\beta\in \Omega},p)\rightarrow (G',\{\prec'_\alpha,\succ'_\alpha,\curlyvee'_{\alpha,\beta}\}_{\alpha,\beta\in \Omega},p')$ of  Hom-NS-family algebras is a linear
map $f:  G \rightarrow  G' $
satisfying $f(x\prec_\alpha y)=f(x)\prec'_\alpha f(y), f(x\succ_\alpha y)=f(x)\succ'_\alpha f(y), f(x\curlyvee_{\alpha,\beta} y)=f(x)\curlyvee'_{\alpha,\beta} f(y)$ for  $x,y\in G$ and $\alpha,\beta\in\Omega$, as well as $f\circ  p =p'\circ f$.

\begin{exam}
(i)  Any  Hom-NS-algebra $(G,\prec,\succ, \curlyvee,p)$  can be considered as a constant  Hom-NS-family algebra $(G,\{\prec_\alpha,\succ_\alpha,\curlyvee_{\alpha,\beta}\}_{\alpha,\beta\in \Omega},p)$,
 where $\prec_\alpha=\prec,\succ_\alpha=\succ $ and $\curlyvee_{\alpha,\beta}=\curlyvee$ for all $\alpha,\beta\in\Omega$.

(ii) Any Hom-dendriform family algebra $(G,\{\prec_\alpha,\succ_\alpha\}_{\alpha\in \Omega},p)$ can be regarded as an
Hom-NS-family algebra $(G,\{\prec_\alpha,\succ_\alpha,\curlyvee_{\alpha,\beta}\}_{\alpha,\beta\in \Omega},p)$ with $\curlyvee_{\alpha,\beta}=0$   for all $\alpha,\beta\in\Omega$.

(iii)
 Let $(G,\{\prec_\alpha,\succ_\alpha,\curlyvee_{\alpha,\beta}\}_{\alpha,\beta\in \Omega})$   be an  NS-family algebra    and $p:G\rightarrow G$ be an NS-family algebra morphism. Define
 $\prec^p_\alpha,\succ^p_\alpha, \curlyvee^p_{\alpha,\beta}:G\otimes G\rightarrow G$ by
$$x\prec^p_\alpha y=p(x)\prec_\alpha p(y),~~~x\succ^p_\alpha y=p(x)\succ_\alpha p(y),~~~x\curlyvee^p_{\alpha,\beta} y=p(x)\curlyvee_{\alpha,\beta} p(y),$$
for all $x,y\in G$ and $\alpha,\beta\in\Omega$. Then $(G,\{\prec^p_\alpha,\succ^p_\alpha, \curlyvee^p_{\alpha,\beta}\}_{\alpha,\beta\in\Omega},p)$   is a  Hom-NS-family algebra, called the
Yau twist of $(G,\{\prec_\alpha,\succ_\alpha,\curlyvee_{\alpha,\beta}\}_{\alpha,\beta\in \Omega})$.
\end{exam}

 \begin{prop} \label{prop:3.6}
 Let   $(G,\{\prec_\alpha,\succ_\alpha,\curlyvee_{\alpha,\beta}\}_{\alpha,\beta\in \Omega},p)$  be a    Hom-NS-family algebra. Then $(G\otimes \mathbb{K}\Omega,\prec,\succ, \curlyvee, \bar{p})$  is a  Hom-NS-algebra, where
 \begin{align*}
 &(a\otimes \alpha)\prec (b\otimes\beta)=(a\prec_\beta b)\otimes\alpha\beta,~~(a\otimes \alpha)\succ (b\otimes\beta)=(a\succ_\alpha b)\otimes\alpha\beta,\\
 &(a\otimes \alpha)\curlyvee (b\otimes\beta)=(a\curlyvee_{\alpha,\beta} b)\otimes\alpha\beta,~~\forall a,b\in G, ~~\alpha,\beta\in \Omega.
  \end{align*}
  Moreover, suppose that $(G',\{\prec'_\alpha,\succ'_\alpha,\curlyvee'_{\alpha,\beta}\}_{\alpha,\beta\in \Omega}, p')$ is another Hom-NS-family algebra. If $f:G\rightarrow G'$ is a morphism of Hom-NS-family algebras,
  then the map
$\bar{f}:G\otimes \mathbb{K}\Omega\rightarrow G'\otimes \mathbb{K}\Omega, a\otimes \alpha\mapsto f(a)\otimes \alpha$ is a morphism between the induced
Hom-NS-algebra.
 \end{prop}

 \begin{proof}
For all $a,b,c\in G$ and $\alpha,\beta,\gamma\in\Omega$, we have
 \begin{align*}
&\bar{p}((a\otimes \alpha)\prec (b\otimes\beta))=\bar{p}((a\prec_\beta b)\otimes\alpha\beta)=p(a\prec_\beta b)\otimes\alpha\beta=p(a)\prec_\beta p(b)\otimes\alpha\beta\\
&=(p(a)\otimes \alpha)\prec (p(b)\otimes\beta)=\bar{p}(a\otimes \alpha)\prec \bar{p}(b\otimes\beta).
 \end{align*}
 Similarly, we get
 \begin{align*}
&\bar{p}((a\otimes \alpha)\succ (b\otimes\beta))=\bar{p}(a\otimes \alpha)\succ \bar{p}(b\otimes\beta),\\
&\bar{p}((a\otimes \alpha)\curlyvee (b\otimes\beta))=\bar{p}(a\otimes \alpha)\curlyvee \bar{p}(b\otimes\beta).
 \end{align*}
Also,  we obtain
  \begin{align*}
&\bar{p}(a\otimes \alpha)\prec \big((b\otimes \beta)\prec (c\otimes \gamma)+(b\otimes \beta)\succ (c\otimes \gamma)+(b\otimes \beta)\curlyvee (c\otimes \gamma)\big)\\
&=(p(a)\otimes \alpha)\prec \big((b\prec_\gamma c+b\succ_\beta c+b\curlyvee_{\beta,\gamma}c)\otimes \beta \gamma\big)\\
&=\big( p(a)\prec_{\beta\gamma}(b\prec_\gamma c+b\succ_\beta c+b\curlyvee_{\beta,\gamma}c)\big)\otimes \alpha\beta \gamma\\
&=\big( (a\prec_{\beta}b)\prec_\gamma p(c) \big)\otimes \alpha\beta \gamma\\
&=((a\prec_{\beta}b)\otimes \alpha\beta) \prec  (p(c)  \otimes  \gamma)\\
&=((a\otimes \alpha)\prec (b\otimes \beta)) \prec \bar{p}(c\otimes \gamma),\\
&((a\otimes \alpha)\succ (b\otimes \beta))\prec \bar{p}(c\otimes \gamma)\\
&=\big((a\succ_\alpha b)\prec_\gamma p(c)\big)\otimes \alpha\beta\gamma=\big(p(a)\succ_\alpha( b\prec_\gamma c)\big)\otimes \alpha\beta\gamma\\
&=\bar{p}(a\otimes \alpha)\succ ((b\otimes \beta)\prec (c\otimes \gamma)),\\
&\big((a\otimes \alpha)\prec (b\otimes \beta)+(a\otimes \alpha)\succ (b\otimes \beta)+(a\otimes \alpha)\curlyvee (b\otimes \beta)\big)\succ \bar{p}(c\otimes \gamma)\\
&=\big((a\prec_\beta b+a\succ_\alpha b+a\curlyvee_{\alpha,\beta} b)\succ_{\alpha\beta} p(c)\big) \otimes \alpha\beta \gamma\\
&=\big(p(a)\succ_\alpha( b \succ_\beta  c)\big) \otimes \alpha\beta \gamma\\
&=\bar{p}(a\otimes \alpha)\succ((b\otimes \beta)\succ (c\otimes \gamma)),\\
&\big((a\otimes \alpha)\prec (b\otimes \beta)+(a\otimes \alpha)\succ (b\otimes \beta)+(a\otimes \alpha)\curlyvee (b\otimes \beta)\big)\curlyvee \bar{p}(c\otimes \gamma)+\\
&\big((a\otimes \alpha)\curlyvee (b\otimes \beta)\big)\prec \bar{p}(c\otimes \gamma) \\
&=\big((a\prec_\beta b+a\succ_\alpha b +a\curlyvee_{\alpha,\beta} b)\curlyvee_{\alpha\beta,\gamma} p(c)\big)\otimes \alpha\beta\gamma+\big((a\curlyvee_{\alpha,\beta} b)\prec_{\gamma} p(c)\big)\otimes \alpha\beta \gamma\\
&=\big(p(a)\succ_{\alpha} (b \curlyvee_{\beta,\gamma} c)\big)\otimes \alpha\beta \gamma+\big((p(a)\curlyvee_{\alpha,\beta\gamma}(b\prec_\gamma c+b\succ_\beta c +b\curlyvee_{\beta,\gamma} c) \big)\otimes \alpha\beta\gamma\\
& =\bar{p}(a\otimes \alpha)\succ\big((b\otimes \beta)\curlyvee (c\otimes \gamma)\big)+\bar{p}(a\otimes \alpha)\curlyvee \big((b\otimes \beta)\prec (c\otimes \gamma)+(b\otimes \beta)\succ (c\otimes \gamma)+(b\otimes \beta)\curlyvee (c\otimes \gamma)\big).
\end{align*}
Hence Eqs. \eqref{3.1}-\eqref{3.5} hold for $(a\otimes \alpha), (b\otimes \beta),(c\otimes \gamma)\in G\otimes \mathbb{K}\Omega$.
That is to say, $(G\otimes \mathbb{K}\Omega,\prec,\succ, \curlyvee, \bar{p})$  is a  Hom-NS-algebra.
The last assertion is obvious.
 \end{proof}

In \cite{Zhang2019}, Zhang and Gao introduced the concept of tridendriform family algebras. Next, we introduce the Hom version of   tridendriform family algebras.
\begin{defn}
  A  Hom-tridendriform family algebra is  a  linear space $G$ together with
a family of bilinear operations  $\{\prec_\alpha,\succ_\alpha:G\times G\rightarrow G\}_{\alpha\in \Omega}$,   a bilinear map $\odot:G\times G\rightarrow G$
and a linear map  $p:G\rightarrow G$ satisfying the following
equations
 \begin{align}
&p(x\prec_\alpha y)= p(x)\prec_\alpha p(y),~~p(x\succ_\alpha y)= p(x)\succ_\alpha p(y),~~p(x\odot y)= p(x)\odot p(y),\label{3.11}\\
& p(x)\prec_{\alpha\beta}(y\prec_\beta z+y\succ_\alpha z+y\odot z)=(x\prec_\alpha y)\prec_\beta p(z),\label{3.12}\\
&(x\succ_\alpha y)\prec_\beta p(z)=p(x)\succ_\alpha(y\prec_\beta z),\label{3.13}\\
&(x\prec_\beta y+x\succ_\alpha y+x\odot y)\succ_{\alpha\beta}p(z)=p(x)\succ_\alpha(y\succ_\beta z),\label{3.14}
 \end{align}
  \begin{align}
&(x\succ_\alpha y)\odot p(z)=p(x)\succ_\alpha (y\odot z), \label{3.15}\\
&(x\prec_\alpha y)\odot p(z)=p(x)\odot(y\succ_\alpha z), \label{3.16}\\
&(x\odot y)\prec_\alpha p(z)=p(x)\odot (y\prec_\alpha z), \label{3.17}\\
&(x\odot y)\odot p(z)=p(x)\odot (y\odot z), \label{3.18}
\end{align}
for all $x,y,z\in G$  and  $\alpha,\beta \in \Omega.$
\end{defn}

\begin{defn}
Let $(L,\mu,p)$ be a  Hom-associative algebra. A collection $\{T_\alpha:L\rightarrow L\}_{\alpha\in\Omega}$ of linear maps is said
to be a Rota-Baxter family operator of weight $\lambda$ if they satisfy
 \begin{align}
&p(T_\alpha x)=T_\alpha p(x), \label{3.19}\\
&T_\alpha x\cdot T_\beta y=T_{\alpha\beta}(T_\alpha x\cdot y+x\cdot T_\beta y+\lambda x\cdot y), \label{3.20}
\end{align}
for all $x,y\in L$  and  $\alpha,\beta\in \Omega.$
\end{defn}
 \begin{prop} \label{prop:3.9}
 Let $(L,\mu,p)$ be a Hom-associative algebra and $\{T_\alpha:L\rightarrow L\}_{\alpha\in\Omega}$  be a    Rota-Baxter family operator of weight $\lambda$.
 Then $(L,\{\prec_\alpha,\succ_\alpha\}_{\alpha\in \Omega},\odot,p)$  is a  Hom-tridendriform family algebra, where
 \begin{align*}
 &x\prec_\alpha y=x\cdot T_\alpha y,~~ x\succ_\alpha y=T_\alpha x\cdot y,~~ x\odot y=\lambda x\cdot y, ~~\forall x,y\in L,~~ \alpha\in \Omega.
  \end{align*}
 \end{prop}

  \begin{proof}
For all $x,y\in L$  and  $\alpha\in \Omega$, we have
  \begin{align*}
&p(x\prec_\alpha y)=p(x\cdot T_\alpha  y)=p(x)\cdot p(T_\alpha  y)=p(x)\cdot T_\alpha p( y)= p(x)\prec_\alpha p(y).
\end{align*}
Similarly, we have
$p(x\succ_\alpha y)= p(x)\succ_\alpha p(y),~~p(x\odot y)= p(x)\odot p(y).$  Next,  for any $x,y,z\in L$  and  $\alpha,\beta\in \Omega$,  we have
\begin{align*}
& p(x)\prec_{\alpha\beta}(y\prec_\beta z+y\succ_\alpha z+y\odot z)=p(x)\cdot T_{\alpha\beta}(y\cdot T_\beta z+T_\alpha y \cdot z+\lambda y\cdot z)\\
&=p(x)\cdot  (T_\alpha y\cdot T_\beta z)=(x\cdot  T_\alpha y)\cdot p(T_\beta z)= (x\cdot  T_\alpha y)\cdot T_\beta p(z)=(x\prec_\alpha y)\prec_\beta p(z), \\
&(x\succ_\alpha y)\prec_\beta p(z)=(T_\alpha x\cdot y)\cdot T_\beta p(z)=(T_\alpha x\cdot y)\cdot p(T_\beta z)=p(T_\alpha x)\cdot (y \cdot  T_\beta z)\\
&=T_\alpha p(x)\cdot (y \cdot  T_\beta z)=p(x)\succ_\alpha(y\prec_\beta z), \\
&(x\prec_\beta y+x\succ_\alpha y+x\odot y)\succ_{\alpha\beta}p(z)=T_{\alpha\beta}(x\cdot T_\beta y+T_\alpha x\cdot y+\lambda x\cdot y)\cdot p(z)\\
&=(T_\alpha x\cdot T_\beta y)\cdot p(z)=p(T_\alpha x)\cdot (T_\beta y \cdot  z)=T_\alpha p(x)\cdot (T_\beta y \cdot  z)=p(x)\succ_\alpha(y\succ_\beta z),\\
&(x\succ_\alpha y)\odot p(z)=\lambda(T_\alpha x\cdot y)\cdot p(z)=\lambda p(T_\alpha x)\cdot (y\cdot  z)= T_\alpha p(x)\cdot \lambda(y\cdot  z)\\
&=p(x)\succ_\alpha (y\odot z),  \\
&(x\prec_\alpha y)\odot p(z)=\lambda(x\cdot T_\alpha y)\cdot p(z)=\lambda p(x)\cdot (T_\alpha y \cdot z)=p(x)\odot(y\succ_\alpha z),  \\
&(x\odot y)\prec_\alpha p(z)=(\lambda x\cdot y)\cdot T_\alpha p(z)=\lambda p(x)\cdot( y\cdot T_\alpha z)=p(x)\odot (y\prec_\alpha z),  \\
&(x\odot y)\odot p(z)=\lambda(\lambda x\cdot y)\cdot p(z)=\lambda p(x)\cdot (\lambda  y\cdot z)=p(x)\odot (y\odot z).
\end{align*}
 Hence Eqs. \eqref{3.12}-\eqref{3.18} hold   and we complete the proof.
  \end{proof}

 \begin{prop} \label{prop:3.10}
 Let $(G,\{\prec_\alpha,\succ_\alpha\}_{\alpha\in \Omega},\odot,p)$  be a  Hom-tridendriform family algebra.
 Then $(G,\{\prec_\alpha,\succ_\alpha, \curlyvee_{\alpha,\beta}\}_{\alpha,\beta\in \Omega}, p)$ is Hom-NS-family algebra, where $\curlyvee_{\alpha,\beta}=\odot$.
 \end{prop}

   \begin{proof}
 Note that the Eqs. \eqref{3.6}-\eqref{3.9} of the Hom-NS-family algebra follow from
the Eqs. \eqref{3.11}-\eqref{3.14} of Hom-tridendriform family algebra. Further, for any $x,y,z\in G$  and  $\alpha,\beta,\gamma\in \Omega$, by  Eqs. \eqref{3.15}-\eqref{3.18},  we have
 \begin{align*}
&p(x)\succ_\alpha(y\curlyvee_{\beta,\gamma}z)+p(x)\curlyvee_{\alpha,\beta\gamma}(y\prec_\gamma z+y\succ_\beta z+y\curlyvee_{\beta,\gamma}z) \\
&=p(x)\succ_\alpha(y\odot z)+p(x)\odot(y\prec_\gamma z)+p(x)\odot(y\succ_\beta z)+p(x)\odot(y\odot z)\\
&=(x \succ_\alpha y)\odot p(z)+(x\odot y)\prec_\gamma p(z)+(x\prec_\beta y)\odot p(z)+(x\odot y)\odot p(z)\\
&=(x\odot y)\prec_\gamma p(z)+(x \succ_\alpha y+x\prec_\beta y+x\odot y)\odot p(z)\\
&=(x\curlyvee_{\alpha,\beta}y)\prec_\gamma p(z)+(x\prec_\beta y+x\succ_\alpha y+x\curlyvee_{\alpha,\beta}y)\curlyvee_{\alpha\beta,\gamma} p(z).
\end{align*}
So Eq. \eqref{3.10} holds and we complete the proof.
  \end{proof}

From  Propositions \ref{prop:3.9} and    \ref{prop:3.10},  the following corollary is clearly established.
 \begin{coro}
 Let $(L,\mu,p)$ be a Hom-associative algebra and $\{T_\alpha:L\rightarrow L\}_{\alpha\in\Omega}$  be a    Rota-Baxter family operator of weight $\lambda$.
 Then $(L,\{\prec_\alpha,\succ_\alpha, \curlyvee_{\alpha,\beta}\}_{\alpha,\beta\in \Omega}, p)$ is Hom-NS-family algebra, where
 \begin{align*}
 &x\prec_\alpha y=x\cdot T_\alpha y,~~ x\succ_\alpha y=T_\alpha x\cdot y,~~ x\curlyvee_{\alpha,\beta} y=\lambda x\cdot y, ~~\forall x,y\in L,~~ \alpha,\beta\in \Omega.
  \end{align*}
 \end{coro}

In the following, we show that an $\Phi$-twisted Rota-Baxter family  operator  induces a  Hom-NS-family
algebra.
  \begin{prop} \label{prop:3.12}
 Let $(L,\mu,p)$ be a Hom-associative algebra, $(V,\mu_l,\mu_r,q)$ be a bimodule.
 and
  $\{R_\alpha:V\rightarrow L\}_{\alpha\in\Omega}$  be an     $\Phi$-twisted Rota-Baxter family  operator.
 Then $(V,\{\prec_\alpha,\succ_\alpha, \curlyvee_{\alpha,\beta}\}_{\alpha,\beta\in \Omega}, q)$ is a Hom-NS-family algebra, where
 \begin{align*}
 &u\prec_\alpha v=u\cdot_r R_\alpha v,~~ u\succ_\alpha v=R_\alpha u\cdot_l v,~~ u\curlyvee_{\alpha,\beta} v=\Phi(R_\alpha u, R_\beta v), ~~\forall u,v\in V,~~ \alpha\in \Omega.
  \end{align*}
  Furthermore, if $\{R_\alpha:V\rightarrow L\}_{\alpha\in\Omega}$ and $\{R'_\alpha:V'\rightarrow L'\}_{\alpha\in\Omega}$ are two  $\Phi$-twisted Rota-Baxter family  operators and $(\phi,\varphi)$ is a morphism between them, then the map $\phi$ is a morphism between the induced Hom-NS-family algebras.
 \end{prop}

\begin{proof}
For all $u,v\in V$ and $\alpha,\beta\in\Omega$, we have
\begin{align*}
&q(u\prec_\alpha v)=q(u\cdot_r R_\alpha v)=q(u)\cdot_r p(R_\alpha v)=q(u)\cdot_rR_\alpha q(v) = q(u)\prec_\alpha q(v),\\
&q(u\succ_\alpha v)=q(R_\alpha u\cdot_l v)=p(R_\alpha u)\cdot_l q(v)=R_\alpha q(u)\cdot_l q(v)= q(u)\succ_\alpha q(v),\\
&q(u\curlyvee_{\alpha,\beta} v)=q(\Phi(R_\alpha u, R_\beta v))=\Phi(p(R_\alpha u), p(R_\beta v))=\Phi(R_\alpha q(u),R_\beta q(v))= q(u)\curlyvee_{\alpha,\beta} q(v).
\end{align*}
 Similarly,  for any $u,v,w\in V$ and $\alpha,\beta,\gamma\in\Omega$, we have
 \begin{align*}
& q(u)\prec_{\alpha\beta}(v\prec_\beta w+v\succ_\alpha w+v\curlyvee_{\alpha,\beta}w)=q(u)\cdot_r R_{\alpha\beta}(v\cdot_r R_\beta w+R_\alpha v\cdot_l w+\Phi(R_\alpha v, R_\beta w))\\
&=q(u)\cdot_r (R_\alpha v \cdot R_\beta w)=(u\cdot_r R_\alpha v) \cdot_r p( R_\beta w)=(u\cdot_r R_\alpha v) \cdot_r R_\beta q(w)\\
&=(u\prec_\alpha v)\prec_\beta q(w), \\
&(u\succ_\alpha v)\prec_\beta q(w)=(R_\alpha u\cdot_l v)\cdot_r R_\beta q(w)=(R_\alpha u\cdot_l v)\cdot_r p(R_\beta w)=p(R_\alpha u)\cdot_l (v\cdot_r R_\beta w)\\
&=R_\alpha q(u)\cdot_l (v\cdot_r R_\beta w)=q(u)\succ_\alpha(v\prec_\beta w), \\
&(u\prec_\beta v+u\succ_\alpha v+u\curlyvee_{\alpha,\beta}v)\succ_{\alpha\beta}q(w)=R_{\alpha\beta}(u\cdot_r R_\beta v+ R_\alpha u\cdot_l v+\Phi(R_\alpha u, R_\beta v))\cdot_l q(w)\\
&=(R_\alpha u\cdot R_\beta v) \cdot_l q(w)=p(R_\alpha u)\cdot_l (R_\beta v  \cdot_l w)=R_\alpha p(u)\cdot_l (R_\beta v  \cdot_l w)=q(u)\succ_\alpha(v\succ_\beta w),  \\
&q(u)\succ_\alpha(v\curlyvee_{\beta,\gamma}w)+q(u)\curlyvee_{\alpha,\beta\gamma}(v\prec_\gamma w+v\succ_\beta w+v\curlyvee_{\beta,\gamma}w) \\
&=R_\alpha q(u)\cdot_l\Phi(R_\beta v, R_\gamma w)+\Phi(R_\alpha q(u),R_{\beta\gamma}(v\cdot_r R_\gamma w+R_\beta v\cdot_l w+\Phi(R_\beta v, R_\gamma w)) )\\
&=p(R_\alpha u)\cdot_l\Phi(R_\beta v, R_\gamma w)+\Phi(p(R_\alpha u), R_\beta v\cdot R_\gamma w)\\
&=\Phi(R_\alpha u,R_\beta v)\cdot_r p(R_\gamma w)+\Phi(R_\alpha u\cdot R_\beta v, p(R_\gamma w))\\
&=\Phi(R_\alpha u,R_\beta v)\cdot_rR_\gamma q(w)+\Phi(R_{\alpha\beta}(u\cdot_r R_\beta v+R_\alpha u\cdot_l v+\Phi(R_\alpha u, R_\beta v)),R_\gamma q(w))\\
&=(u\curlyvee_{\alpha,\beta}v)\prec_\gamma q(w)+(u\prec_\beta v+u\succ_\alpha v+u\curlyvee_{\alpha,\beta}v)\curlyvee_{\alpha\beta,\gamma} q(w).
 \end{align*}
Thus, $(V,\{\prec_\alpha,\succ_\alpha, \curlyvee_{\alpha,\beta}\}_{\alpha,\beta\in \Omega}, q)$ is a Hom-NS-family algebra.
 \end{proof}

 The following corollaries can be obtained   from Proposition \ref{prop:3.12}.

  \begin{coro} \label{coro:3.13}
 Let $(L,\mu,p)$ be a Hom-associative algebra,  $(V,\mu_l,\mu_r,q)$ be a bimodule and $\Phi\in C^2_{HA}(L,V)$  be a 2-cocycle.  If $R:V\rightarrow L$  is an     $\Phi$-twisted Rota-Baxter    operator,
 then $(V, \prec,\succ, \curlyvee, q)$ is a Hom-NS-algebra, where
 \begin{align*}
 &u\prec v=u\cdot_r R v,~~ u\succ v=Ru\cdot_l v,~~ u\curlyvee v=\Phi(R u, R v), ~~\forall u,v\in V.
  \end{align*}
 \end{coro}

  \begin{coro}
 Let $(L,\mu,p)$ be a Hom-associative algebra and $\{N_\alpha:L\rightarrow L\}_{\alpha\in\Omega}$  be a    Nijenhuis family operator .
 Then $(L,\{\prec_\alpha,\succ_\alpha, \curlyvee_{\alpha,\beta}\}_{\alpha,\beta\in \Omega}, p)$ is a Hom-NS-family algebra, where
 \begin{align*}
 &x\prec_\alpha y=x\cdot  N_\alpha y,~~ x\succ_\alpha y=N_\alpha x\cdot  y,~~ x\curlyvee_{\alpha,\beta} y=-N_{\alpha\beta} (x\cdot y), ~~\forall x,y\in L,~~ \alpha,\beta\in \Omega.
  \end{align*}
 \end{coro}
   \begin{proof}
 By Example \ref{exam:2.8},  for any $ x,y\in L,~~ \alpha,\beta\in \Omega$, we have
   \begin{align*}
&x\prec_\alpha y=x\cdot N_\alpha y=x\cdot_r^N  (y\otimes\alpha)=x\cdot_r^N  Id_\alpha(y),\\
 &x\succ_\alpha y=N_\alpha x\cdot  y= (x\otimes\alpha)\cdot_l^N y= Id_\alpha(x)\cdot_l^N y\\
 &  \text{and}\\
 & x\curlyvee_{\alpha,\beta} y=-N_{\alpha\beta} (x\cdot y)=\Phi_N(x\otimes\alpha,y\otimes\beta)=\Phi_N(Id_\alpha(x),Id_\beta(y)),
  \end{align*}
where $\{Id_\alpha :L\rightarrow L\otimes\mathbb{K}\Omega\}$ is an $\Phi_N$-twisted Rota-Baxter family operator.  According to Proposition \ref{prop:3.12}, $(L,\{\prec_\alpha,\succ_\alpha, \curlyvee_{\alpha,\beta}\}_{\alpha,\beta\in \Omega}, p)$ is a Hom-NS-family algebra.
  \end{proof}

 Let $(L,\mu, p)$  be another  Hom-associative algebra and $(V,\mu_l,\mu_r,q)$ be a bimodule of
$L$.
Suppose that  $\{R_\alpha:V\rightarrow L\}_{\alpha\in\Omega}$  is an     $\Phi$-twisted Rota-Baxter family  operator, then  by Proposition \ref{prop:3.12},
$(V,\{\prec_\alpha,\succ_\alpha, \curlyvee_{\alpha,\beta}\}_{\alpha,\beta\in \Omega}, q)$ is a Hom-NS-family algebra.
Hence, by Proposition \ref{prop:3.6},   $(V\otimes \mathbb{K}\Omega,\prec,\succ, \curlyvee, \bar{p})$  is a  Hom-NS-algebra, where
 \begin{align}
\left\{ \begin{array}{lll}
 (u\otimes \alpha)\prec (v\otimes\beta)=(u\prec_\beta v)\otimes\alpha\beta=(u\cdot_r R_\alpha v)\otimes\alpha\beta,\\
 (u\otimes \alpha)\succ (v\otimes\beta)=(u\succ_\alpha v)\otimes\alpha\beta=(R_\alpha u\cdot_l v)\otimes\alpha\beta,\\
 (u\otimes \alpha)\curlyvee (v\otimes\beta)=(u\curlyvee_{\alpha,\beta} v)\otimes\alpha\beta=\Phi(R_\alpha u, R_\beta v)\otimes\alpha\beta.
\end{array}  \right. \label{3.21}
 \end{align}

 On the other hand, by Proposition \ref{prop:2.9}, the $\Phi$-twisted Rota-Baxter family  operator $\{R_\alpha:V\rightarrow L\}_{\alpha\in\Omega}$ induces an  $\bar{\Phi}$-twisted Rota-Baxter   operator $\bar{R}:V\otimes \mathbb{K}\Omega\rightarrow L\otimes \mathbb{K}\Omega, R(u\otimes \alpha)=R_\alpha u\otimes\alpha$. Thus, we have a  Hom-NS-algebra structure on $V\otimes \mathbb{K}\Omega$ by
Corollary \ref{coro:3.13}, which coincides with the one given by \eqref{3.21}. As a summary,
 we have shown that $\Phi$-twisted Rota-Baxter family  operators,  Hom-NS-family algebras, $\bar{\Phi}$-twisted Rota-Baxter   operators and  Hom-NS-algebras are
closely related in the sense of commutative diagram of categories as follows:
 $$\aligned
\xymatrix{
  \text{$\Phi$-twisted Rota-Baxter family  operators} \ar[r]  \ar[d] & \text{Hom-NS-family algebras} \ar[d]\\
   \text{$\bar{\Phi}$-twisted Rota-Baxter   operators} \ar[r]& \text{Hom-NS-algebras.}  }
 \endaligned$$

\section{Cohomology    of  twisted Rota-Baxter family  operators} \label{sec:Cohomology}
\def\theequation{\arabic{section}.\arabic{equation}}
\setcounter{equation} {0}

 In \cite{Aguiar2020},  Aguiar introduced the concept of associative algebra relative to a semigroup ($\Omega$-associative algebra). Next, we introduce the Hom version of  associative algebras relative to a semigroup.
\begin{defn}
 A Hom-$\Omega$-associative algebra (Hom-associative algebra relative to a semigroup $\Omega$) is a vector space $G$ together with a collection of bilinear operations $\{\cdot_{\alpha,\beta}:G\otimes G\rightarrow G\}_{\alpha,\beta\in \Omega}$ and a  linear map  $  p:G\rightarrow G$ satisfying
  \begin{align*}
 &  p(x\cdot_{\alpha,\beta}y)=p(x)\cdot_{\alpha,\beta}p(y),\\
 &p(x)\cdot_{\alpha,\beta\gamma}(y\cdot_{\beta,\gamma}z)=(x\cdot_{\alpha,\beta}y)\cdot_{\alpha\beta,\gamma} p(z), ~~\forall ~x,y,z\in G, \alpha,\beta,\gamma\in \Omega.
 \end{align*}
\end{defn}

   \begin{prop} \label{prop:4.2}
 Let  $(G,\{\prec_\alpha,\succ_\alpha, \curlyvee_{\alpha,\beta}\}_{\alpha,\beta\in \Omega}, p)$ be a Hom-NS-family algebra.
 Then $(G,\{*_{\alpha,\beta}\}_{\alpha,\beta\in \Omega},p)$ is a Hom-$\Omega$-associative algebra, where
  \begin{align*}
 x*_{\alpha,\beta}y=x\prec_\beta y+x\succ_\alpha y+x \curlyvee_{\alpha,\beta}y, ~~\forall ~x,y\in G, \alpha,\beta\in \Omega.
 \end{align*}
 \end{prop}
    \begin{proof}
 By Eqs. \eqref{3.6},    we have $p (x*_{\alpha,\beta}y)=p(x)*_{\alpha,\beta}p(y)$. Further,  for all ~$x,y,z\in G, \alpha,\beta,\gamma\in \Omega$, by Eqs. \eqref{3.7}-\eqref{3.10},  we have
   \begin{align*}
&p(x)*_{\alpha,\beta\gamma}(y*_{\beta,\gamma}z)\\
&=p(x)*_{\alpha,\beta\gamma}(y\prec_\gamma z+y\succ_\beta z+y \curlyvee_{\beta,\gamma}z)\\
&=p(x)\prec_{\beta\gamma}(y\prec_\gamma z+y\succ_\beta z+y \curlyvee_{\beta,\gamma}z)+p(x)\succ_{\alpha}(y\prec_\gamma z+y\succ_\beta z+y \curlyvee_{\beta,\gamma}z)+\\
&~~~~~p(x)\curlyvee_{\alpha,\beta\gamma}(y\prec_\gamma z+y\succ_\beta z+y \curlyvee_{\beta,\gamma}z)\\
&=(x\prec_{\beta}y)\prec_\gamma p(z)+p(x)\succ_{\alpha}(y\prec_\gamma z)+p(x)\succ_{\alpha}(y\succ_\beta z)+p(x)\succ_{\alpha}(y \curlyvee_{\beta,\gamma}z)+\\
&~~~~~p(x)\curlyvee_{\alpha,\beta\gamma}(y\prec_\gamma z+y\succ_\beta z+y \curlyvee_{\beta,\gamma}z)\\
&=(x\prec_{\beta}y)\prec_\gamma p(z)+(x\succ_{\alpha} y)\prec_\gamma p(z)+(x\prec_\beta y+x\succ_\alpha y+x \curlyvee_{\alpha,\beta}y)\succ_{\alpha\beta}p(z)+\\
&~~~~ (x \curlyvee_{\alpha,\beta}y)\prec_{\gamma}p(z)+(x\prec_{\beta} y+x\succ_\alpha y+x \curlyvee_{\alpha,\beta}y)\curlyvee_{\alpha\beta,\gamma}p(z)\\
&=(x*_{\alpha,\beta}y)\prec_\gamma p(z) +(x*_{\alpha,\beta}y)\succ_{\alpha\beta}p(z)+ (x*_{\alpha,\beta}y)\curlyvee_{\alpha\beta,\gamma}p(z)\\
&=(x*_{\alpha,\beta}y)*_{\alpha\beta,\gamma} p(z).
  \end{align*}
Therefore, $(G,\{*_{\alpha,\beta}\}_{\alpha,\beta\in \Omega},p)$ is a Hom-$\Omega$-associative algebra.
  \end{proof}

  Combining Propositions \ref{prop:3.12} and \ref{prop:4.2}, we get the following result.

  \begin{coro} \label{coro:4.3}
  Let $(L,\mu, p)$  be a   Hom-associative algebra, $(V,\mu_l,\mu_r,q)$ be a bimodule of
$L$  and  $\{R_\alpha:V\rightarrow L\}_{\alpha\in\Omega}$  be an     $\Phi$-twisted Rota-Baxter family  operator.
 Then $(V,\{*_{\alpha,\beta}\}_{\alpha,\beta\in \Omega},q)$ is a Hom-$\Omega$-associative algebra, where
  \begin{align*}
   u*_{\alpha,\beta} v= R_\alpha u \cdot_l v+u\cdot_rR_\beta v+\Phi(R_\alpha u, R_\beta v), ~~\forall ~u,v\in V, \alpha,\beta\in \Omega.
 \end{align*}
  \end{coro}

\begin{defn}
Let  $(G,\{\cdot_{\alpha,\beta}\}_{\alpha,\beta\in \Omega},p)$ be a Hom-$\Omega$-associative algebra. A bimodule over it consists of a
Hom-vector space $(V,q)$ together with a collection
 \begin{align*}
 \left\{\begin{array}{cccc}
 \mu_{l_{\alpha,\beta}}:G\otimes V\rightarrow V,~~~ (x,u)\mapsto x\cdot_{l_{\alpha,\beta}}u,\\
 \mu_{r_{\alpha,\beta}}:V \otimes G \rightarrow V, ~~~(u,x)\mapsto u\cdot_{r_{\alpha,\beta}}x
 \end{array}\right\}_{\alpha,\beta\in\Omega}
 \end{align*}
of bilinear maps satisfying, for $x,y\in G, u\in V$ and $\alpha,\beta,\gamma\in \Omega$,
  \begin{align}
 &  q(x\cdot_{l_{\alpha,\beta}}u)=p(x)\cdot_{l_{\alpha,\beta}}q(u),~~q(u\cdot_{r_{\alpha,\beta}}x)=q(u)\cdot_{l_{\alpha,\beta}}p(x),\label{4.1}\\
 &q(u)\cdot_{r_{\alpha,\beta\gamma}}(x\cdot_{\beta,\gamma}y)=(u\cdot_{r_{\alpha,\beta}}x)\cdot_{r_{\alpha\beta,\gamma}} p(y), \label{4.2}\\ &p(x)\cdot_{l_{\alpha,\beta\gamma}}(u\cdot_{r_{\beta,\gamma}}y)=(x\cdot_{l_{\alpha,\beta}}u)\cdot_{r_{\alpha\beta,\gamma}} p(y),\label{4.3}\\
 &p(x)\cdot_{l_{\alpha,\beta\gamma}}(y\cdot_{l_{\beta,\gamma}}u)=(x\cdot_{\alpha,\beta}y)\cdot_{l_{\alpha\beta,\gamma}} q(u). \label{4.4}
 \end{align}
 \end{defn}

Next,  we assume that $\Omega$ is a semigroup with unit $1\in \Omega$. The unital condition of $\Omega$ is only useful in the coboundary
operator of the cohomology at the degree 0 level.
   \begin{lemma} \label{lemma:4.4}
 Let  $(G,\{\cdot_{\alpha,\beta}\}_{\alpha,\beta\in \Omega},p)$ be a Hom-$\Omega$-associative algebra and $(V,\mu_{l_{\alpha,\beta}},\mu_{r_{\alpha,\beta}}q)$ be a bimodule over it.
 For each $n\geq 0$, we define $C^n_{\Omega}(G,V)$ by
   \begin{align*}
 &C^0_{\Omega}(G,V)=\{u\in V~|~q(u)=u\},\\
 &C^{n\geq 1}_{\Omega}(G,V)=\{f=\{f_{\alpha_1,\ldots,\alpha_n}\in \mathrm{Hom}(G^{\otimes n},V)\}_{\alpha_1,\ldots,\alpha_n\in\Omega}~|~q\circ f=f\circ p^{\otimes n}\}.
 \end{align*}
  The  map $\delta_\Omega^n: C^{n}_{\Omega}(G,V)\rightarrow C^{n+1}_{\Omega}(G,V)$  given by
  \begin{align*}
  (\delta_\Omega^0(u))_\alpha(x)=&x\cdot_{l_{\alpha,1}} u-u\cdot_{r_{1,\alpha}}x,\\
  (\delta_\Omega^n(f))_{\alpha_1,\ldots,\alpha_{n+1}}(x_1,\ldots,x_{n+1})=&p^{n-1}(x_1)\cdot_{l_{\alpha_1,\alpha_2\cdots\alpha_{n+1}}}f_{\alpha_2,\ldots,\alpha_{n+1}}(x_2,\ldots,x_{n+1})+\\
  &(-1)^{n+1}f_{\alpha_1,\ldots,\alpha_{n}}(x_1,\ldots,x_{n})\cdot_{r_{\alpha_1\cdots\alpha_n,\alpha_{n+1}}}p^{n-1}(x_{n+1})+&\\
 & \sum_{i=1}^n(-1)^i f_{\alpha_1,\ldots,\alpha_i\alpha_{i+1},\ldots,\alpha_{n+1}}(p(x_1),\ldots,x_i\cdot_{\alpha_i,\alpha_{i+1}}x_{i+1},\ldots, p(x_{n+1})).&
 \end{align*}
 for any $u\in C^0_{\Omega}(G,V)$ and $f\in C^{n}_{\Omega}(G,V)$. Then $\delta_\Omega^{n+1} \circ \delta_\Omega^{n}=0.$ That is, for $n\geq 0$, $(C^{n\geq 0}_{\Omega}(G,V),\delta_\Omega^n)$ is a cochain complex. The corresponding cohomology groups are called the cohomology of the Hom-$\Omega$-associative algebra
$G$ with coefficients in $V$.
 \end{lemma}

   \begin{theorem} \label{theorem:4.5}
 Let $(L,\mu, p)$  be a   Hom-associative algebra, $(V,\mu_l,\mu_r,q)$ be a bimodule of
$L$  and  $\{R_\alpha:V\rightarrow L\}_{\alpha\in\Omega}$  be an     $\Phi$-twisted Rota-Baxter family  operator.
Define a collection of
bilinear maps
  \begin{align*}
  \left\{\begin{array}{cc}
 \lhd_{\alpha,\beta}:V\otimes L\rightarrow L,~~ u\lhd_{\alpha,\beta}x= R_\alpha u\cdot x-R_{\alpha\beta}(u\cdot_r x)-R_{\alpha\beta}\Phi(R_\alpha u, x), \\
 \rhd_{\alpha,\beta}:L\otimes V\rightarrow L,~~x\rhd_{\alpha,\beta}u= x \cdot R_\beta u-R_{\alpha\beta}(x\cdot_l u)-R_{\alpha\beta}\Phi(x, R_\beta u)
 \end{array}\right\}_{\alpha,\beta\in\Omega}
 \end{align*}
 Then $(L,\{\lhd_{\alpha,\beta},\rhd_{\alpha,\beta}\}_{\alpha,\beta\in \Omega},p)$ is a bimodule of the Hom-$\Omega$-associative algebra $(V,\{*_{\alpha,\beta}\}_{\alpha,\beta\in \Omega},q)$
 introduced in   Corollary  \ref{coro:4.3}.
 \end{theorem}

    \begin{proof}
For all ~$x\in L, u\in V$ and $  \alpha,\beta \in \Omega$,     we have
   \begin{align*}
&p(u\lhd_{\alpha,\beta}x)=p(R_\alpha u\cdot x-R_{\alpha\beta}(u\cdot_r x)-R_{\alpha\beta}\Phi(R_\alpha u, x))\\
&=R_\alpha q(u)\cdot p(x)-R_{\alpha\beta}(q(u)\cdot_r p(x))-R_{\alpha\beta}\Phi(R_\alpha q(u), p(x))\\
&=q(u)\lhd_{\alpha,\beta}p(x).
  \end{align*}
  Similarly,  $p(x\rhd_{\alpha,\beta}u)=p(x)\rhd_{\alpha,\beta}q(u)$. Hence \eqref{4.1} holds. Now we verify \eqref{4.2}-\eqref{4.4}. For any ~$x\in L, u,v\in V$ and $  \alpha,\beta,\gamma \in \Omega$,   by Eqs. \eqref{2.2}
  -- \eqref{2.4},    we have
 \begin{small}
  \begin{align*}
 &p(x)\rhd_{\alpha,\beta\gamma}(u*_{\beta,\gamma}v)=p(x)\rhd_{\alpha,\beta\gamma}(R_\beta u \cdot_l v+u\cdot_rR_\gamma v+\Phi(R_\beta u, R_\gamma v))\\
 &=p(x) \cdot R_{\beta\gamma} (R_\beta u \cdot_l v+u\cdot_rR_\gamma v+\Phi(R_\beta u, R_\gamma v))-R_{\alpha\beta\gamma}(p(x)\cdot_l (R_\beta u \cdot_l v+u\cdot_rR_\gamma v+\Phi(R_\beta u, R_\gamma v)))-\\
 &~~~~R_{\alpha\beta\gamma}\Phi(p(x), R_{\beta\gamma} (R_\beta u \cdot_l v+u\cdot_rR_\gamma v+\Phi(R_\beta u, R_\gamma v)))\\
 &=p(x) \cdot  (R_\beta u\cdot R_\gamma v) -R_{\alpha\beta\gamma}\big(p(x)\cdot_l (R_\beta u \cdot_l v+u\cdot_rR_\gamma v+\Phi(R_\beta u, R_\gamma v))+\Phi(p(x),  R_\beta u\cdot R_\gamma v)\big) \\
  &=(x \cdot  R_\beta u)\cdot p(R_\gamma v) -R_{\alpha\beta\gamma}\big((x\cdot R_\beta u) \cdot_l q(v)+(x\cdot_l u)\cdot_r p(R_\gamma v)+\Phi(x,R_\beta u)\cdot_rp(R_\gamma v) +\Phi(x\cdot R_\beta u, p(R_\gamma v))\big) \\
  &=(x \cdot R_\beta u)\cdot R_\gamma q(v)-R_{\alpha\beta\gamma}((x \cdot R_\beta u)\cdot_l q(v))-R_{\alpha\beta}(x\cdot_l u+\Phi(x, R_\beta u)) \cdot R_\gamma q(v)+\\
  &~~R_{\alpha\beta\gamma}\big(R_{\alpha\beta}(x\cdot_l u+\Phi(x, R_\beta u))\cdot_l q(v)+\Phi(R_{\alpha\beta}(x\cdot_l u+\Phi(x, R_\beta u)), R_\gamma q(v))\big)-R_{\alpha\beta\gamma}\Phi(x \cdot R_\beta u, R_\gamma q(v))\\
   &=(x \cdot R_\beta u-R_{\alpha\beta}(x\cdot_l u)-R_{\alpha\beta}\Phi(x, R_\beta u)) \cdot R_\gamma q(v)-R_{\alpha\beta\gamma}((x \cdot R_\beta u-R_{\alpha\beta}(x\cdot_l u)-R_{\alpha\beta}\Phi(x, R_\beta u))\cdot_l q(v))-\\
  &~~~R_{\alpha\beta\gamma}\Phi((x \cdot R_\beta u-R_{\alpha\beta}(x\cdot_l u)-R_{\alpha\beta}\Phi(x, R_\beta u)), R_\gamma q(v))\\
  &=(x \cdot R_\beta u-R_{\alpha\beta}(x\cdot_l u)-R_{\alpha\beta}\Phi(x, R_\beta u))\rhd_{\alpha\beta,\gamma} q(v) \\
 &=(x\rhd_{\alpha,\beta}u)\rhd_{\alpha\beta,\gamma} q(v),
  \end{align*}
  \begin{align*}
 &q(u)\lhd_{\alpha,\beta\gamma}(x\rhd_{\beta,\gamma}v)=q(u)\lhd_{\alpha,\beta\gamma}(x \cdot R_\gamma v-R_{\beta\gamma}(x\cdot_l v)-R_{\beta\gamma}\Phi(x, R_\gamma v))\\
 &=R_\alpha q(u)\cdot (x \cdot R_\gamma v-R_{\beta\gamma}(x\cdot_l v)-R_{\beta\gamma}\Phi(x, R_\gamma v))-R_{\alpha\beta\gamma}(q(u)\cdot_r (x \cdot R_\gamma v-R_{\beta\gamma}(x\cdot_l v)-R_{\beta\gamma}\Phi(x, R_\gamma v)))-\\
 &~~~R_{\alpha\beta\gamma}\Phi(R_\alpha q(u), x \cdot R_\gamma v-R_{\beta\gamma}(x\cdot_l v)-R_{\beta\gamma}\Phi(x, R_\gamma v))\\
  &=(R_\alpha u\cdot x) \cdot R_\gamma q(v)-R_{\alpha\beta\gamma}\big(R_\alpha q(u)\cdot_l(x\cdot_l v+\Phi(x, R_\gamma v))+q(u)\cdot_rR_{\beta\gamma}(x\cdot_l v+\Phi(x, R_\gamma v))+\\
  &~~~\Phi(R_\alpha q(u),  R_{\beta\gamma}(x\cdot_l v+\Phi(x, R_\gamma v)))\big)-R_{\alpha\beta\gamma}(q(u)\cdot_r (x \cdot R_\gamma v-R_{\beta\gamma}(x\cdot_l v)-R_{\beta\gamma}\Phi(x, R_\gamma v)))-\\
 &~~~R_{\alpha\beta\gamma}\Phi(R_\alpha q(u), x \cdot R_\gamma v-R_{\beta\gamma}(x\cdot_l v)-R_{\beta\gamma}\Phi(x, R_\gamma v))\\
   &=(R_\alpha u\cdot x) \cdot R_\gamma q(v)-R_{\alpha\beta\gamma}\big(R_\alpha q(u)\cdot_l(x\cdot_l v+\Phi(x, R_\gamma v))+q(u)\cdot_r (x \cdot R_\gamma v)+\Phi(R_\alpha q(u), x \cdot R_\gamma v)\big)\\
 &=(R_\alpha u\cdot x) \cdot R_\gamma q(v)-R_{\alpha\beta\gamma}\big((R_\alpha u \cdot  x)\cdot_l q(v)+p(R_\alpha u)\cdot_l\Phi(x, R_\gamma v)+ (u \cdot_r  x) \cdot_r R_\gamma q(v)+\Phi(p(R_\alpha u), x \cdot R_\gamma v)\big)\\
 &=(R_\alpha u\cdot x) \cdot R_\gamma q(v)-R_{\alpha\beta\gamma}\big((R_\alpha u \cdot  x)\cdot_l q(v)+ (u \cdot_r  x) \cdot_r R_\gamma q(v)+\Phi(R_\alpha u,x)\cdot_r p(R_\gamma v)+\Phi(R_\alpha u\cdot x ,p(R_\gamma v))\big)\\
&=(R_\alpha u\cdot x)\cdot R_\gamma  q(v)-R_{\alpha\beta\gamma}((R_\alpha u \cdot  x)\cdot_l q(v))-R_{\alpha\beta}(u\cdot_r x+\Phi(R_\alpha u, x))\cdot R_\gamma  q(v)+\\
&~~~R_{\alpha\beta\gamma} \big(R_{\alpha\beta}(u\cdot_r x+\Phi(R_\alpha u, x))\cdot_l  q(v)+\Phi(R_{\alpha\beta}(u\cdot_r x+\Phi(R_\alpha u, x)), R_\gamma  q(v))\big)-R_{\alpha\beta\gamma}\Phi(R_\alpha u\cdot x, R_\gamma  q(v))\\
&=(R_\alpha u\cdot x-R_{\alpha\beta}(u\cdot_r x)-R_{\alpha\beta}\Phi(R_\alpha u, x))\cdot R_\gamma  q(v)-R_{\alpha\beta\gamma}((R_\alpha u\cdot x-R_{\alpha\beta}(u\cdot_r x)-R_{\alpha\beta}\Phi(R_\alpha u, x))\cdot_l  q(v))-\\
&~~~R_{\alpha\beta\gamma}\Phi((R_\alpha u\cdot x-R_{\alpha\beta}(u\cdot_r x)-R_{\alpha\beta}\Phi(R_\alpha u, x)), R_\gamma  q(v))\\
 &=(R_\alpha u\cdot x-R_{\alpha\beta}(u\cdot_r x)-R_{\alpha\beta}\Phi(R_\alpha u, x))\rhd_{\alpha\beta,\gamma} q(v)\\
 &=(u\lhd_{\alpha,\beta}x)\rhd_{\alpha\beta,\gamma} q(v)
   \end{align*}
   and
  \begin{align*}
 &q(u)\lhd_{\alpha,\beta\gamma}(v\lhd_{\beta,\gamma}x)=q(u)\lhd_{\alpha,\beta\gamma}(R_\beta v\cdot x-R_{\beta\gamma}(v\cdot_r x)-R_{\beta\gamma}\Phi(R_\beta v, x))\\
 &=R_\alpha q(u)\cdot  (R_\beta v\cdot x-R_{\beta\gamma}(v\cdot_r x)-R_{\beta\gamma}\Phi(R_\beta v, x))-R_{\alpha\beta\gamma}(q(u)\cdot_r  (R_\beta v\cdot x-R_{\beta\gamma}(v\cdot_r x)-R_{\beta\gamma}\Phi(R_\beta v, x)))-\\
 &~~~R_{\alpha\beta\gamma}\Phi(R_\alpha q(u),  R_\beta v\cdot x-R_{\beta\gamma}(v\cdot_r x)-R_{\beta\gamma}\Phi(R_\beta v, x))\\
  &=R_\alpha q(u)\cdot  (R_\beta v\cdot x)-  R_{\alpha\beta\gamma}\big(R_\alpha q(u)\cdot_l(v\cdot_r x+\Phi(R_\beta v, x))+ q(u)\cdot_r R_{\beta\gamma}(v\cdot_r x+\Phi(R_\beta v, x))+\\
  &~~~\Phi(R_\alpha q(u), R_{\beta\gamma}(v\cdot_r x+\Phi(R_\beta v, x)))\big)-R_{\alpha\beta\gamma}(q(u)\cdot_r  (R_\beta v\cdot x-R_{\beta\gamma}(v\cdot_r x)-R_{\beta\gamma}\Phi(R_\beta v, x)))-\\
 &~~~R_{\alpha\beta\gamma}\Phi(R_\alpha q(u),  R_\beta v\cdot x-R_{\beta\gamma}(v\cdot_r x)-R_{\beta\gamma}\Phi(R_\beta v, x))\\
   &=R_\alpha q(u)\cdot  (R_\beta v\cdot x)-  R_{\alpha\beta\gamma}\big(R_\alpha q(u)\cdot_l(v\cdot_r x+\Phi(R_\beta v, x))+q(u)\cdot_r  (R_\beta v\cdot x)+\Phi(R_\alpha q(u),  R_\beta v\cdot x)\big)\\
 &=(R_\alpha u \cdot   R_\beta v)\cdot p(x)-  R_{\alpha\beta\gamma}\big((R_\alpha u \cdot_l v)\cdot_r p(x)+p(R_\alpha u)\cdot_l\Phi(R_\beta v, x) +(u\cdot_r  R_\beta v)\cdot p(x)+\Phi(p(R_\alpha u),  R_\beta v\cdot x)\big)\\
  &=(R_\alpha u \cdot   R_\beta v)\cdot p(x)-R_{\alpha\beta\gamma}\big(( R_\alpha u \cdot_l v)\cdot_r p(x) +(u\cdot_rR_\beta v)\cdot_r p(x)+\Phi(R_\alpha u, R_\beta v)\cdot_r p(x)+\Phi(R_\alpha u \cdot   R_\beta v, p(x))\big)\\
 &=R_{\alpha\beta} ( R_\alpha u \cdot_l v+u\cdot_rR_\beta v+\Phi(R_\alpha u, R_\beta v))\cdot p(x)-R_{\alpha\beta\gamma}(( R_\alpha u \cdot_l v+u\cdot_rR_\beta v+\Phi(R_\alpha u, R_\beta v))\cdot_r p(x))-\\
 &~~~R_{\alpha\beta\gamma}\Phi(R_{\alpha\beta} ( R_\alpha u \cdot_l v+u\cdot_rR_\beta v+\Phi(R_\alpha u, R_\beta v)), p(x))\\
 &=( R_\alpha u \cdot_l v+u\cdot_rR_\beta v+\Phi(R_\alpha u, R_\beta v))\lhd_{\alpha\beta,\gamma} p(x)\\
 &=(u*_{\alpha,\beta}v)\lhd_{\alpha\beta,\gamma} p(x).
 \end{align*}
   \end{small}
This completes the proof.
  \end{proof}

Finally, we consider the cochain complex of the Hom-$\Omega$-associative algebra $(V,\{*_{\alpha,\beta}\}_{\alpha,\beta\in \Omega},q)$ with coefficients in the bimodule $(L,\{\lhd_{\alpha,\beta},\rhd_{\alpha,\beta}\}_{\alpha,\beta\in \Omega},p)$
given in the above theorem. More precisely, for all $n\geq 0$, we define
  \begin{align*}
 &C^0_{R}(V,L)=\{x\in L~|~p(x)=x\},\\
 &C^{n\geq 1}_{R}(V,L)=\{f=\{f_{\alpha_1,\ldots,\alpha_n}\in \mathrm{Hom}(V^{\otimes n},L)\}_{\alpha_1,\ldots,\alpha_n\in\Omega}~|~p\circ f=f\circ q^{\otimes n}\}
 \end{align*}
and the differential by $\delta^n_{R}: C^{n}_{R}(V,L)\rightarrow C^{n+ 1}_{R}(V,L)$ by
  \begin{align}
  &(\delta^0_{R}(x))_\alpha(u)\nonumber\\
  =&R_\alpha u \cdot x-R_\alpha(u\cdot_rx)-R_\alpha\Phi(R_\alpha u,x)-x\cdot R_\alpha u+R_\alpha(x\cdot_l u)+R_\alpha\Phi(x,R_\alpha u),\label{4.5}\\
 & (\delta^n_{R}(f))_{\alpha_1,\ldots,\alpha_{n+1}}(u_1,\ldots, u_{n+1})\nonumber\\
 =&R_{\alpha_1}q^{n-1}(u_1)\cdot f_{\alpha_2,\ldots,\alpha_{n+1}}(u_2,\ldots, u_{n+1})-R_{\alpha_1\ldots\alpha_{n+1}}(q^{n-1}(u_1)\cdot_r f_{\alpha_2,\ldots,\alpha_{n+1}}(u_2,\ldots, u_{n+1}))-\nonumber\\
 &R_{\alpha_1\ldots\alpha_{n+1}}\Phi(R_{\alpha_1}q^{n-1}(u_1), f_{\alpha_2,\ldots,\alpha_{n+1}}(u_2,\ldots, u_{n+1}))+\nonumber\\
 &\sum_{i=1}^n(-1)^if_{\alpha_1,\ldots,\alpha_{i}\alpha_{i+1},\ldots,\alpha_{n+1}}\Big(q(u_1),\ldots,q(u_{i-1}),R_{\alpha_i} u_i\cdot_l   u_{i+1}+ u_i\cdot_r R_{\alpha_{i+1}} u_{i+1}+\nonumber\\
 &\Phi(R_{\alpha_i} u_i, R_{\alpha_{i+1}} u_{i+1}),q(u_{i+1}),\ldots,  q(u_{n+1})\Big)+\nonumber\\
 &(-1)^{n+1} \Big(f_{\alpha_1,\ldots,\alpha_{n}}(u_1,\ldots, u_{n})\cdot R_{\alpha_{n+1}}q^{n-1}(u_{n+1})-R_{\alpha_{1}\cdots\alpha_{n+1}} \big(f_{\alpha_1,\ldots,\alpha_{n}}(u_1,\ldots, u_{n})\cdot_l q^{n-1}(u_{n+1})\big)  -\nonumber\\
 &  R_{\alpha_{1}\cdots\alpha_{n+1}}\Phi(f_{\alpha_1,\ldots,\alpha_{n}}(u_1,\ldots, u_{n}), R_{\alpha_{n+1}}q^{n-1}(u_{n+1})      \Big), \label{4.6}
 \end{align}
for all $x\in C^0_{R}(V,L)$ and $f\in C^{n\geq 1}_{R}(V,L)$. Let $Z^n_{R}(V,L)=\{f\in C^{n}_{R}(V,L)~|~\delta^n_{R}(f)=0\}$ be the space of $n$-cocycles and
$B^n_{R}(V,L)=\{\delta^{n-1}_{R}(f)~|~f\in C^{n-1}_{R}(V,L)\}$ the space of $n$-coboundaries, for $n\geq 0$. The
corresponding cohomology groups
$H^n_{R}(V,L)=\frac{Z^n_{R}(V,L)}{B^n_{R}(V,L)}$
are called the cohomology of the     $\Phi$-twisted Rota-Baxter family  operator $\{R_\alpha:V\rightarrow L\}_{\alpha\in\Omega}$.

   \begin{remark}
 (i)  Let $\{R_\alpha:V\rightarrow L\}_{\alpha\in\Omega}$  be a relative  Rota-Baxter family  operator (see Example \ref{exam:2.5}). Then $(V,\{*_{\alpha,\beta}\}_{\alpha,\beta\in \Omega},q)$ is a Hom-$\Omega$-associative algebra, where
  \begin{align*}
   u*_{\alpha,\beta} v= R_\alpha u \cdot_l v+u\cdot_rR_\beta v, ~~\forall ~u,v\in V, \alpha,\beta\in \Omega.
 \end{align*}
Moreover, $(L,\{\lhd_{\alpha,\beta},\rhd_{\alpha,\beta}\}_{\alpha,\beta\in \Omega},p)$ is a bimodule of   $(V,\{*_{\alpha,\beta}\}_{\alpha,\beta\in \Omega},q)$, where
  \begin{align*}
\left\{\begin{array}{cc}
\lhd_{\alpha,\beta}:V\otimes L\rightarrow L,~~ u\lhd_{\alpha,\beta}x= R_\alpha u\cdot x-R_{\alpha\beta}(u\cdot_r x), \\
\rhd_{\alpha,\beta}:L\otimes V\rightarrow L,~~x\rhd_{\alpha,\beta}u= x \cdot R_\beta u-R_{\alpha\beta}(x\cdot_l u)
\end{array} \right\}_{\alpha,\beta\in\Omega}
 \end{align*}
 The corresponding cohomology is called the cohomology of the relative  Rota-Baxter family  operator $\{R_\alpha:V\rightarrow L\}_{\alpha\in\Omega}$.

 (ii) The cohomology of the $\Phi$-twisted Rota-Baxter family  operator defined as above is the
cohomology of the $\Phi$-twisted Rota-Baxter   operator   when the semigroup $\Omega$ is trivial.
 \end{remark}

    \begin{remark}
 The cohomology theory for $\Phi$-twisted Rota-Baxter family  operators enjoys certain
functorial properties. Let $\{R_\alpha:V\rightarrow L\}_{\alpha\in\Omega}$    and  $\{R'_\alpha:V\rightarrow L\}_{\alpha\in\Omega}$ be   two  $\Phi$-twisted Rota-Baxter family operators on a bimodule $(V,\mu_l,\mu_r,q)$ over the Hom-associative algebra  $L$, and $(\phi,\varphi)$    a morphism from $\{R_\alpha:V\rightarrow L\}_{\alpha\in\Omega}$ to $\{R'_\alpha:V\rightarrow L\}_{\alpha\in\Omega}$ in which $\phi$ is invertible.
 Let $C^{n}_{R}(V,L)$ denote the space
of $n$-cochains of $\Phi$-twisted Rota-Baxter family  operators. Define a map $\xi^n:C^{n}_{R}(V,L)\rightarrow C^{n}_{R'}(V,L) $ by
 \begin{align*}
 (\xi^n(f))_{\alpha_1,\ldots,\alpha_{n+1}}(v_1,\ldots,v_{n+1})=\varphi(f_{\alpha_1,\ldots,\alpha_{n+1}}(\phi^{-1}(v_1),\ldots,\phi^{-1}(v_{n+1})),
 \end{align*}
 where $  f\in C^{n}_{R}(V,L), \alpha_1,\ldots,\alpha_{n+1}\in \Omega $ and $ v_1,\ldots,v_{n+1}\in V.$ Then it is straightforward
to deduce that $\xi^\bullet$ is a cochain map from the   complex $(\oplus_{n=0}^{\infty}C^{n}_{R}(V,L),\delta^\bullet_{R})$
to the   complex $(\oplus_{n=0}^{\infty}C^{n}_{R'}(V,L),\delta^\bullet_{R'})$, that is,
$\xi^{n+1} \circ \delta^{n}_{R}=\delta^{n}_{R'} \circ \xi^n.$
In other words, the following diagram is commutative:
$$\aligned
\xymatrix{
  C^{n}_{R}(V,L)\ar[r]^-{\delta^n_{R}}\ar[d]^-{\xi^n}& C^{n+1}_{R}(V,L)\ar[d]^{\xi^{n+1}}\\
  C^{n}_{R'}(V,L)\ar[r]^-{\delta^n_{R'}}& C^{n+1}_{R'}(V,L).}
 \endaligned$$
 Consequently, $\xi^n$ induces a  morphism $\xi^n_{\ast}$ from the cohomology group
 $H^{n}_{R}(V,L)$ to $H^{n}_{R'}(V,L)$.
 \end{remark}

 \section{Deformations  of  twisted Rota-Baxter family  operators} \label{sec:Deformations}
\def\theequation{\arabic{section}.\arabic{equation}}
\setcounter{equation} {0}

In this section, we study the    deformations of $\Phi$-twisted Rota-Baxter family  operators on Hom-associative algebras  following the classical approaches proposed  by Gerstenhaber \cite{Gerstenhaber63,Gerstenhaber64}.
We show that the cohomology of  $\Phi$-twisted Rota-Baxter family  operators introduced above govern such    deformations.

\begin{defn}
 Let $\mathcal{R}=\{R_\alpha:V\rightarrow L\}_{\alpha\in\Omega}$     be  an  $\Phi$-twisted Rota-Baxter family operator on a bimodule $(V,\mu_l,\mu_r,q)$ over the Hom-associative algebra  $(L,\mu, p)$.
 An infinitesimal   deformation of $\mathcal{R}$ consists of a sum $\mathcal{R}^t=\mathcal{R}+t\mathcal{R}^1$, where $\mathcal{R}^1=\{R^1_\alpha:V\rightarrow L\}_{\alpha\in\Omega}$, such that $\mathcal{R}^t$ is an $\Phi$-twisted Rota-Baxter family operator  on   a bimodule $(V,\mu_l,\mu_r,q)$ over   $(L,\mu, p)$, for all $t$. In such a case, we say that $\mathcal{R}^1$ generates an infinitesimal or a
linear deformation of $\mathcal{R}$.
\end{defn}

It follows that $\mathcal{R}^t=\mathcal{R}+t\mathcal{R}^1$ is an infinitesimal   deformation of $\mathcal{R}$   if and only if for any $u,v\in V$ and $\alpha,\beta\in \Omega$,
\begin{align*}
&R^t_\alpha q(u)=p( R^t_\alpha u), \\
&R^t_\alpha u\cdot R^t_\beta v=R^t_{\alpha\beta}(R^t_\alpha u \cdot_l v+u\cdot_rR^t_\beta v+\Phi(R^t_\alpha u, R^t_\beta v)).
\end{align*}
By equating coefficients of $t$  from both sides, we have
\begin{align}
&R^1_\alpha q(u)=p( R^1_\alpha u),\label{5.1} \\
&R^1_\alpha u\cdot R_\beta v+R_\alpha u\cdot R^1_\beta v=R^1_{\alpha\beta}\big(R_\alpha u \cdot_l v+u\cdot_rR_\beta v+\Phi(R_\alpha u, R_\beta v)\big)+\nonumber\\
&~~~R_{\alpha\beta}\big(R^1_\alpha u \cdot_l v+u\cdot_rR^1_\beta v+\Phi(R^1_\alpha u, R_\beta v)+\Phi(R_\alpha u, R^1_\beta v)\big).\label{5.2}
\end{align}

It follows from Eqs.\eqref{5.1} and  \eqref{5.2} that $\mathcal{R}^1=\{R^1_\alpha:V\rightarrow L\}_{\alpha\in\Omega}\in C^{ 1}_{R}(V,L)$ and $(\delta_{R}^1(\mathcal{R}^1))_{\alpha,\beta}(u,v)=0$,   respectively.
 That is $\mathcal{R}^1 $ is a 1-cocycle in the cochain complex of the $\Phi$-twisted Rota-Baxter family operator  $\{R_\alpha:V\rightarrow L\}_{\alpha\in\Omega}$.
The 1-cocycle $\mathcal{R}^1$ is called the infinitesimal of the deformation $\mathcal{R}^t$.

Given a Hom-NS-family algebra  $(G,\{\prec_\alpha,\succ_\alpha,\curlyvee_{\alpha,\beta}\}_{\alpha,\beta\in \Omega}, p)$. Let $\{\prec^1_\alpha,\succ^1_\alpha,\curlyvee^1_{\alpha,\beta}\}_{\alpha,\beta\in \Omega}$ be a family
of bilinear operations. If for any $t$, the multiplications $\{\prec^t_\alpha,\succ^t_\alpha,\curlyvee^t_{\alpha,\beta}\}_{\alpha,\beta\in \Omega}$ defined by
$$\Big\{\prec^t_\alpha=\prec_\alpha+t\prec^1_\alpha,~~\succ^t_\alpha=\succ_\alpha+t\succ^1_\alpha,~~\curlyvee^t_{\alpha,\beta}=\curlyvee_{\alpha,\beta}+t\curlyvee^1_{\alpha,\beta}\Big\}_{\alpha,\beta\in\Omega}$$
give a  Hom-NS-family algebra structure on $(G,p)$, we say that the pair $\{\prec^1_\alpha,\succ^1_\alpha,\curlyvee^1_{\alpha,\beta}\}_{\alpha,\beta\in \Omega}$ generates
a   infinitesimal deformation of the Hom-NS-family algebra  $(G,\{\prec_\alpha,\succ_\alpha,\curlyvee_{\alpha,\beta}\}_{\alpha,\beta\in \Omega}, p)$.

In the following, we show that an infinitesimal deformation of an $\Phi$-twisted Rota-Baxter family operator induces an
infinitesimal deformation of the corresponding Hom-NS-family algebra  structure on the module.

  \begin{prop}
  If $\mathcal{R}^1=\{R^1_\alpha:V\rightarrow L\}_{\alpha\in\Omega}$ generates an infinitesimal deformation of an $\Phi$-twisted Rota-Baxter family operator $\{R_\alpha:V\rightarrow L\}_{\alpha\in\Omega}$ on a bimodule $(V,\mu_l,\mu_r,q)$ over the Hom-associative algebra  $(L,\mu, p)$, then the multiplications
   \begin{align*}
\left\{\begin{array}{cccc}
 u\prec^t_\alpha v=u\cdot_r R_\alpha v+tu\cdot_r R^1_\alpha v,\\
   u\succ^t_\alpha v=R_\alpha u\cdot_l v+tR^1_\alpha u\cdot_l v,\\
   u\curlyvee^t_{\alpha,\beta} v=\Phi(R_\alpha u, R_\beta v)+t\big(\Phi(R^1_\alpha u, R_\beta v)+\Phi(R_\alpha u, R^1_\beta v)\big)
  \end{array} \right\}_{\alpha,\beta\in\Omega}
  \end{align*}
  define an infinitesimal deformation of the corresponding Hom-NS-family algebra structure on
$(V,q)$.
 \end{prop}

  \begin{coro}
  If $\mathcal{R}^1=\{R^1_\alpha:V\rightarrow L\}_{\alpha\in\Omega}$ generates an infinitesimal deformation of an $\Phi$-twisted Rota-Baxter family operator $\{R_\alpha:V\rightarrow L\}_{\alpha\in\Omega}$ on a bimodule $(V,\mu_l,\mu_r,q)$ over the Hom-associative algebra  $(L,\mu, p)$, then the multiplication
    \begin{align*}
 \Big\{  u*^t_{\alpha,\beta} v=u*_{\alpha,\beta} v+ t\big(R^1_\alpha u \cdot_l v+u\cdot_rR^1_\beta v+\Phi(R^1_\alpha u, R_\beta v)+\Phi(R_\alpha u, R^1_\beta v)\big)\Big \}_{\alpha,\beta\in\Omega}
 \end{align*}
 defines an infinitesimal deformation of the corresponding Hom-$\Omega$-associative algebra structure on
$(V,q)$.
  \end{coro}

\begin{defn}
 Let $\mathcal{R}=\{R_\alpha:V\rightarrow L\}_{\alpha\in\Omega}$     be  an  $\Phi$-twisted Rota-Baxter family operator on a bimodule $(V,\mu_l,\mu_r,q)$ over the Hom-associative algebra  $(L,\mu, p)$.
 Two infinitesimal deformations $\mathcal{R}^t=\mathcal{R}+t\mathcal{R}^1$ and $ \bar{\mathcal{R}}^t=\mathcal{R}+t \bar{\mathcal{R}}^1$
 of $\mathcal{R}$ are equivalent if there exists an element $x\in L$ such that $p(x)=x$ and the pair
 $\big(\phi^t=Id_V+t(\mu_{l_x}-\mu_{r_x}+\Phi(x,R_\alpha-)-\Phi(R_\alpha-,x)), \varphi^t=Id_L+t(\mu(x,-)-\mu(-,x))\big)$ is
a morphism from  $\mathcal{R}^t$  to $  \bar{\mathcal{R}}^t$.
\end{defn}

Let us recall from Definition \ref{defn:2.4} that the pair $ (\phi^t,\varphi^t)$ is
a morphism of $\Phi$-twisted Rota-Baxter family operators from  $\mathcal{R}^t$  to $ \bar{\mathcal{R}}^t$ if the following conditions are satisfied

(i) The map $\varphi^t$ is a Hom-associative algebra homomorphism.

(ii) $\varphi^t\circ \mathcal{R}^t=\bar{\mathcal{R}}^t\circ\phi^t$.


(iii) $\phi^t\circ q=q\circ \phi^t$.

(iv) $\phi^t\circ\Phi=\Phi\circ(\varphi^t\otimes\varphi^t)$.

(v) $\phi^t\circ\mu_l=\mu_l \circ(\varphi^t\otimes\phi^t).$

(vi) $\phi^t\circ \mu_r=\mu_r \circ(\phi^t\otimes\varphi^t)$.

\begin{theorem}
  Let $\mathcal{R}^t=\mathcal{R}+t\mathcal{R}^1$ and $ \bar{\mathcal{R}}^t=\mathcal{R}+t \bar{\mathcal{R}}^1$ be two equivalent infinitesimal
deformations of an $\Phi$-twisted Rota-Baxter family operator $\mathcal{R}=\{R_\alpha:V\rightarrow L\}_{\alpha\in\Omega}$. Then $\mathcal{R}^1=\{R^1_\alpha:V\rightarrow L\}_{\alpha\in\Omega}$ and $\bar{\mathcal{R}}^1=\{\bar{R}^1_\alpha:V\rightarrow L\}_{\alpha\in\Omega}$ defines the same cohomology class in $H^1_R(V,L)$.
\end{theorem}

\begin{proof}
Assuming that $\big(\phi^t=Id_V+t(\mu_{l_x}-\mu_{r_x}+\Phi(x,R_\alpha-)-\Phi(R_\alpha-,x)), \varphi^t=Id_L+t(\mu(x,-)-\mu(-,x))\big)$ is
a morphism from  $\mathcal{R}^t$  to $  \bar{\mathcal{R}}^t$,  from the above condition (ii), for any $u\in V$,  we   obtain
  \begin{align*}
 &R_\alpha u +tR^1_\alpha u +tx\cdot R_\alpha u +t^2x\cdot R^1_\alpha u-tR_\alpha u \cdot x-t^2R^1_\alpha u \cdot x\\
 =&R_\alpha u +t\big(R_\alpha (x\cdot_l u)-R_\alpha (u\cdot_r x)+R_\alpha \Phi(x,R_\alpha u)-R_\alpha \Phi(R_\alpha u, x) \big) +\\
 &t\bar{R}^1_\alpha u  +t^2\big(\bar{R}^1_\alpha (x\cdot_l u)-\bar{R}^1_\alpha (u\cdot_r x)+\bar{R}^1_\alpha \Phi(x,R_\alpha u)+\bar{R}^1_\alpha \Phi(R_\alpha u, x) \big).
 \end{align*}
 Comparing coefficients of $t^1$ on both sides of above equation,  we have
  \begin{align*}
 &R^1_\alpha u +x\cdot R_\alpha u -R_\alpha u \cdot x=R_\alpha (x\cdot_l u)-R_\alpha (u\cdot_r x)+R_\alpha \Phi(x,R_\alpha u)-R_\alpha \Phi(R_\alpha u, x)+\bar{R}^1_\alpha u .
 \end{align*}
 Further by  Eq. \eqref{4.5}, we   get $R^1_\alpha u-\bar{R}^1_\alpha u=(\delta_R^0(x))_\alpha(u)$, which implies that $\mathcal{R}^1$ and $\bar{\mathcal{R}}^1$ belong to the same cohomology class in $H^1_R(V,L)$.
\end{proof}

Regarding the formal deformation of $\Phi$-twisted Rota-Baxter family operators, analogous results are obtained.
 \begin{remark}
 Let  $\mathcal{R}^t=\mathcal{R}+\sum_{i=1}^{\infty}\mathcal{R}^it^i$ be any formal deformation of   the $\Phi$-twisted Rota-Baxter family operator $\mathcal{R}=\{R_\alpha:V\rightarrow L\}_{\alpha\in\Omega}$. Then the infinitesimal $\mathcal{R}^1$ is a 1-cocycle in the
cochain complex of $\mathcal{R}$.  Moreover, the corresponding cohomology class depends
only on the equivalence class of the deformation $\mathcal{R}^t$.
 \end{remark}

\begin{defn}
Let $\mathcal{R}=\{R_\alpha:V\rightarrow L\}_{\alpha\in\Omega}$     be  an  $\Phi$-twisted Rota-Baxter family operator on a bimodule $(V,\mu_l,\mu_r,q)$ over the Hom-associative algebra  $(L,\mu, p)$.
An infinitesimal deformation $\mathcal{R}^t=\mathcal{R}+t\mathcal{R}^1$  is said
to be trivial if it is equivalent to the deformation $\bar{\mathcal{R}}^0=\mathcal{R}$.
\end{defn}

Let $\big(\phi^t=Id_V+t(\mu_{l_x}-\mu_{r_x}+\Phi(x,R_\alpha-)-\Phi(R_\alpha-,x)), \varphi^t=Id_L+t(\mu(x,-)-\mu(-,x))\big)$ be
a morphism from  $\mathcal{R}^t$  to $  \bar{\mathcal{R}}^0=\mathcal{R}$.
  Then   $\varphi^t$ is a Hom-associative algebra homomorphism  of $L$, i.e., the following equations hold:
  \begin{align}
 &(x\cdot a)\cdot (x\cdot b)  -(x\cdot a)\cdot (b\cdot x)-(a\cdot x)\cdot (x\cdot b)+(a\cdot x)\cdot (b\cdot x)=0,~~~ \forall  a,b\in L. \label{5.3}
 \end{align}
By  $\phi^t \Phi(a,b)=\Phi(\varphi^t(a), \varphi^t(b))$,  we have
\begin{align}
\left\{ \begin{array}{lll}
 x\cdot_l  \Phi(a,b)- \Phi(a,b)\cdot_r x+ \Phi(x, R_\alpha\Phi(a,b))- \Phi(R_\alpha\Phi(a,b),x)  \\
 =\Phi(x\cdot a-a\cdot x,b)+\Phi(a,x\cdot b-b\cdot x),  \\
 \Phi(x\cdot a-a\cdot x, x\cdot b-b\cdot x)=0.
 \end{array}  \right. \label{5.4}
 \end{align}
Finally, by $\phi^t(a\cdot_l u)= \varphi^t(a)\cdot_l\phi^t(u)$ and $\phi^t(u\cdot_r a)= \phi^t(u)\cdot_r\varphi^t(a)$,  we have
\begin{align}
\left\{ \begin{array}{lll}
 x\cdot_l  (a\cdot_l u)-(a\cdot_l u)\cdot_r x+  \Phi(x, R_\alpha (a\cdot_l u))-\Phi(R_\alpha (a\cdot_l u),x) \\
 =(x\cdot a-a\cdot x)\cdot_lu+a\cdot_l(x\cdot_lu-u\cdot_r x+\Phi(x, R_\alpha u)-\Phi(R_\alpha u, x)),  \\
 (x\cdot a-a\cdot x)\cdot_l(x\cdot_lu-u\cdot_l x+\Phi(x, R_\alpha u)-\Phi(R_\alpha u, x))=0,
 \end{array}  \right. \label{5.5}
 \end{align}
 \begin{align}
 \left\{ \begin{array}{lll}
  x\cdot_l  (u\cdot_r a)-(u\cdot_r a)\cdot_r x+  \Phi(x, R_\alpha (u\cdot_r a))-\Phi(R_\alpha (u\cdot_r a),x) \\
 =u \cdot_r(x\cdot a-a\cdot x)+(x\cdot_lu-u\cdot_r x+\Phi(x, R_\alpha u)-\Phi(R_\alpha u, x)) \cdot_r a, \\
 (x\cdot_lu-u\cdot_l x+\Phi(x, R_\alpha u)-\Phi(R_\alpha u, x))\cdot_r (x\cdot a-a\cdot x)=0.
\end{array}  \right. \label{5.6}
 \end{align}

Next we define Nijenhuis elements associated to an $\Phi$-twisted Rota-Baxter family operator $\mathcal{R}=\{R_\alpha:V\rightarrow L\}_{\alpha\in\Omega}$   in a way
that a trivial deformation of $\mathcal{R}$ induces a Nijenhuis element.

\begin{defn}
An element $x\in L$ is called a Nijenhuis element associated
to  $\mathcal{R}$  if $x$  such that $p(x)=x$,
  \begin{align*}
 & x\cdot(u\lhd_{\alpha,\beta}x-x\rhd_{\alpha,\beta}u)-(u\lhd_{\alpha,\beta}x-x\rhd_{\alpha,\beta}u)\cdot x=0
 \end{align*}
 and the Eqs. \eqref{5.3}--\eqref{5.6} hold. The set of all Nijenhuis elements associated with $\mathcal{R}$ is denoted by $\mathrm{Nij}(\mathcal{R})$.
\end{defn}

Finally,  we give a sufficient condition for the rigidity of an  $\Phi$-twisted Rota-Baxter family operator in terms of Nijenhuis elements.

\begin{theorem}
 Let $\mathcal{R}=\{R_\alpha:V\rightarrow L\}_{\alpha\in\Omega}$     be  an  $\Phi$-twisted Rota-Baxter family operator on a bimodule $(V,\mu_l,\mu_r,q)$ over the Hom-associative algebra  $(L,\mu, p)$.
  If $Z_R^1(V,L)=\delta_R^0(\mathrm{Nij}(\mathcal{R}))$, then $\mathcal{R}$ is
rigid.
\end{theorem}

\begin{proof}
Let $\mathcal{R}^t=\mathcal{R}+\sum_{i=1}^{\infty}\mathcal{R}^it^i$ be any formal deformation of $\mathcal{R}$. According to the above discussion, the linear term $\mathcal{R}^1$
is a 1-cocycle in the cohomology of $\mathcal{R}$, i.e., $\mathcal{R}^1\in Z_R^1(V,L)$. Thus, by the hypothesis, there is a Nijenhuis
element $x\in \mathrm{Nij}(\mathcal{R})$ such that $\mathcal{R}^1=\delta_R^0(x)$. We take
$$\phi^t=Id_V+t(\mu_{l_x}-\mu_{r_x}+\Phi(x,R_\alpha-)-\Phi(R_\alpha-,x))~~\text{and}~~ \varphi^t=Id_L+t(\mu(x,-)-\mu(-,x)) $$
and define $\bar{\mathcal{R}}^t= \varphi^t\circ \mathcal{R}^t\circ (\phi^t)^{-1}$.  Then $\bar{\mathcal{R}}^t$ is a formal deformation equivalent to ${\mathcal{R}}^t$. It follows that
  \begin{align*}
 & \bar{\mathcal{R}}^t(u)\\
 &=(Id_L+t(\mu(x,-)-\mu(-,x)))\circ\\
 &~~~(R_\alpha+tR_\alpha^1+\text{power of} ~t^{\geq 2})(u-t(x\cdot_l u-u\cdot_rx+\Phi(x,R_\alpha u)-\Phi(R_\alpha u,x)))\\
 &=(Id_L+t(\mu(x,-)-\mu(-,x)))\big(R_\alpha u+tR_\alpha^1u-tR_\alpha (x\cdot_l u-u\cdot_rx+\Phi(x,R_\alpha u)-\Phi(R_\alpha u,x))+\\
  &~~~\text{power of} ~t^{\geq 2}\big)\\
  &= R_\alpha u+t(R_\alpha^1u+x\cdot R_\alpha u-R_\alpha u\cdot x-R_\alpha (x\cdot_l u-u\cdot_rx+\Phi(x,R_\alpha u)-\Phi(R_\alpha u,x)))+\\
  &~~~\text{power of} ~t^{\geq 2} \\
  &= R_\alpha u+t^2 \bar{{R}}^2_\alpha u+\cdots,
 \end{align*}
where the last equation follows by ${R}^1_\alpha=(\delta_R^0(x))_\alpha$. This indicates that the coefficient
of $t$ in the expression of $\bar{\mathcal{R}}^t$ is trivial. By repeating this process,  we ascertain that $\mathcal{R}^t$ is
equivalent to $\mathcal{R}$,  indicating that $\mathcal{R}$ is indeed rigid.
\end{proof}

{{\bf Acknowledgment.}  \\
The paper is  supported by the NSF of China (No. 12161013) and Guizhou Provincial Basic Research Program (Natural
Science) (No. ZK[2023]025).

\end{document}